\documentclass[journal]{IEEEtran}

\usepackage[cmex10]{amsmath}
\usepackage{amsfonts,mathtools}
\usepackage{algorithmic}
\usepackage[]{algorithm}
\usepackage{array}
\usepackage{url}
\usepackage{bm}
\usepackage{graphicx}
\usepackage{mathrsfs}
\usepackage{amsthm} 
\usepackage{amsmath,bm}
\usepackage{subfigure}
\usepackage{hyperref}
\usepackage{cite}
\usepackage{color}
\usepackage{float}
\usepackage{caption}
\usepackage{epstopdf}
\theoremstyle{remark}
\newtheorem{rem}{\bf Remark}
\theoremstyle{plain}
\newtheorem{lem}{Lemma}
\newtheorem{theorem}{Theorem}
 \theoremstyle{plain}

\newtheorem{assumption}{Assumption}

\def \z {{\mathbf z}}
\def \x {{\mathbf x}}
\def \y {{\mathbf y}}
\def \X {{\mathbf X}}

\def \N {{\mathcal N}}
\def \d {{\mathbb R}}
\def \set  {{\mathcal S}}

\def \Rn {{\mathbb{R}}}

\providecommand{\abs}[1]{\lvert#1\rvert}                      
\newcommand{\ind}{1\hspace{-1.6mm}1}                          
\newcommand{\norm}[1]{\ensuremath{\left\|#1\right\|}}          

\def \cG {{\mathcal{G}}}
\def \cK {{\mathcal{K}}}
\def \cN {{\mathcal{N}}}
\def \cE {{\mathcal{E}}}
\def \cX {{\mathcal{X}}}
\def \pj {{\mathrm{prox}_j}}
\def \pk {{\mathrm{prox}_k}}
\def \bt {{\boldsymbol{\theta}}}

\def \X {{\mathbf{X}}}
\def \x {\mathbf{x}}

\def \z {\mathbf{z}}

\def \S {\mathcal{S}}

\begin{document}

		\setlength{\abovedisplayskip}{1pt}
		\setlength{\belowdisplayskip}{1pt}
 	\captionsetup[figure]{skip=-.5pt}

\title{Asynchronous Optimization Over Heterogeneous Networks via Consensus ADMM}

\author{ Sandeep Kumar, \IEEEmembership{ Student member, IEEE}, Rahul Jain, and Ketan Rajawat, \IEEEmembership{Member, IEEE} 
\thanks{Sandeep Kumar and Ketan Rajawat are with the Dept. of Electrical Engineering, IIT Kanpur, Kanpur, UP-208016, India; emails: \texttt{\{sandkr,ketan\}@iitk.ac.in}. Rahul jain is with Numerify.com, Bangalore Karnatka-560025, India; email \texttt{rahuljain2104@gmail.com}}}

\maketitle
\begin{abstract}
This paper considers the distributed optimization of a sum of locally observable, non-convex functions. The optimization is performed over a multi-agent networked system, and each local function depends only on a subset of the variables. An asynchronous and distributed alternating directions method of multipliers (ADMM) method that allows the nodes to defer or skip the computation and transmission of updates is proposed in the paper. The proposed algorithm utilizes different approximations in the update step, resulting in proximal and majorized ADMM variants. Both variants are shown to converge to a local minimum, under certain regularity conditions. The proposed asynchronous algorithms are also applied to the problem of cooperative localization in wireless ad hoc networks, where it is shown to outperform the other state-of-the-art localization algorithms. 
\end{abstract}

\IEEEpeerreviewmaketitle

\section{Introduction}
Multi-agent networked systems arise in a number of engineering disciplines such as tactical ad hoc networks\cite{martin2005distributed}, environmental monitoring networks \cite{speranzon2006distributed}, multi-robot control and tracking \cite{localization_example, hero}, internet-scale monitoring\cite{mateos2013dynamic,cortes2009distributed}, and large-scale learning \cite{forero2008consensus}. The estimation, resource allocation, and network control tasks required in these applications are often formulated as distributed optimization problems, where each node is associated with a local cost function, determined from set of local and possibly private measurements \cite{admm_defn,duchi2012dual}. The nodes must nevertheless cooperate in order to minimize the network objective function, which is the sum of local costs. The algorithm design becomes challenging in the absence of a centralized controller or a fusion center, where nodal interactions are limited to their neighborhoods, and global network state information is largely unavailable \cite{parikh2013proximal, admm_defn}. In high-dimensional problems, even local message passing may be prohibitive, since each node may only be interested in a subset of the optimization variables. 

In general, distributed optimization algorithms are designed either in the primal \cite{duchi2012dual,distributedSubgradient,lobel2011distributed}, \cite[Chap. 10]{palomar2010convex} or dual domain \cite{admm_defn,terelius2011decentralized}. A popular dual approach is the distributed alternating direction method of multipliers (ADMM), where the nodal variables are decoupled through the introduction of the consensus constraints, and the updates are carried out in the dual domain \cite{admm_defn,order}. High-dimensional problems are handled through the so-called \emph{general-form consensus} formulation, where local updates depend only on a subset of optimization variables\cite[Chap. 7]{admm_defn}. The ADMM algorithm is also applicable to a class of non-convex problems, where it has been shown to converge to a local minimum \cite{admm,hong2014convergence}.

Practical networks are also heterogeneous with respect to their processing powers, energy availability, and communication capabilities, giving rise to asynchrony  \cite{parallel,tsitsiklis1986distributed}. Indeed, a key feature required of the distributed algorithms is their tolerance to processing and communication delays arising due to slow or energy-starved nodes \cite{li2013distributed,agarwal2011distributed,liu2013asynchronous,scutari2008asynchronous,nedic2001distributed,order,async_admm,srivastava2011distributed,zhao2015asynchronous}. In general, distributed optimization algorithms such as ADMM must be appropriately modified to allow updates to be skipped or delayed, and the convergence of the asynchronous variants is neither obvious, nor guaranteed \cite[Chap. 6]{parallel,admm}. Asynchronous or randomized variants of the distributed ADMM are well-known for the case when the cost functions are convex; see e.g., \cite{order,async_admm,chang2014proximal} and references therein. Likewise, if a master node or a fusion center is available, the asynchronous distributed ADMM variant proposed in \cite{admm}, is applicable to non-convex cost functions and has been shown to converge to a local optimum. 

This work considers the non-convex general-form consensus optimization problem arising in a multi-agent networked system. The first contribution is the development of an asynchronous and distributed ADMM framework that runs without a fusion center. Two variants, namely, the proximal and the majorized ADMM, are proposed, each handling the non-convex objective function in a different way. While the ADMM iterations at all nodes still occur according to a common schedule, both algorithms allow the nodes to, at times, skip the computationally intensive steps and/or the transmission of updates. As the second contribution, it is shown that both variants converge to a local minimum, under certain regularity conditions. The convergence analysis reveals that with appropriately chosen parameters, the algorithms can tolerate any bounded level of asynchrony.  Finally, the third contribution is the application of the proposed ADMM algorithm to the problem of cooperative localization for distributed networks \cite{localization_example, hero}. Detailed simulations and comparisons with existing distributed and asynchronous localization algorithms are carried out, establishing the superior performance of the ADMM algorithm. 

This paper is organized as follows. The problem formulation and related examples are provided in Sec. \ref{probfor}. The proposed proximal asynchronous ADMM and the associated convergence results are presented in Sec. \ref{asec1}, while the proposed majorized asynchronous ADMM is detailed in Sec. \ref{madmms}. Finally, the simulation results for the proposed algorithm applied to the cooperative localization problem are provided in Sec. \ref{sim} and sec. \ref{conclusion} concludes the paper.

\section{Problem Formulation}\label{probfor}
{ This section details the partially separable non-convex problem formulation considered here, and motivates the need for a distributed optimization algorithm via two examples. Before describing the problem at hand, some background regarding the general multi-agent optimization problem is first presented.}

\subsection{Background}
Consider a network represented by the undirected graph $\cG = (\cK,\cE)$, where $\cK : = \{1, 2, \ldots, K\}$ denotes the set of agents or nodes and $\cE$, the set of edges that represent  communication links. A node $k \in \cK$ may only communicate with its neighbors $\N'_k := \{j | (j,k) \in \cE\}$. Consider first the general multi-agent problem where the nodes want to cooperatively solve the following optimization problem:
\begin{align}\label{prob1}   
\min   f( \x)  := & \sum _{ k=1 }^{ K }{ g_{ k }( \x)  }  + h(\x) \quad
\text{ s. t. }   \x\in  \cX 
\end{align}
where $\x \in \d^{N \times 1}$, $g_k: \d^{N} \rightarrow \d$ for $k = 1$, $\ldots$, $K$ are differentiable, possibly non-convex functions, and $h: \d^N \rightarrow \d$ is a convex but not necessarily differentiable function. The set $\cX$ is closed, convex, and compact. Within the distributed setting considered here, the function $g_k(\x)$ is local to node $k$, and the network has no central coordinator or fusion center. 

In general, non-convex problems such as \eqref{prob1} are solved in a distributed manner using the first order gradient or subgradient descent, dual methods such as ADMM, convex relaxation (such as semidefinite relaxation), successive convex optimization, or leveraging the problem structure; see \cite{nonconvex_examples} and references therein. These approaches result in algorithms that are parallelizable to various extents, with different computational and message passing requirements. Of particular interest here is the high-dimensional regime, where large $N$ prohibits nodes from operating over and exchanging the full vector $\x$. { To this end, the next section considers a partially separable form of \eqref{prob1} which is amenable to a distributed optimization algorithm.}

\begin{figure}[tbh]
	\centering
	\includegraphics[scale=.55]{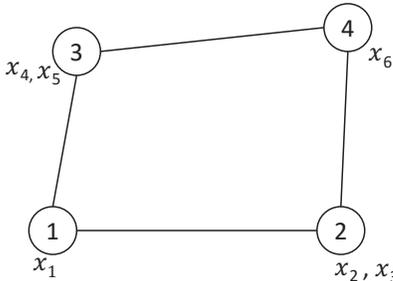}
	\caption{An example of a four-node network.}
	\label{network}
\end{figure}

\begin{figure}[tbh]
	\centering
	\includegraphics[scale=.55]{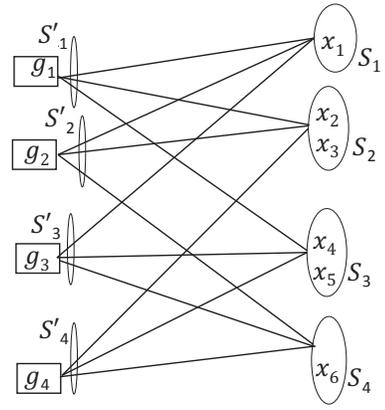}
	\caption{Factor graph representation for the objective function of \eqref{prob2}.}
	\label{separable:fig}
\end{figure}

\subsection{Partially separable form}\label{separable}
Consider a special case of \eqref{prob1} where the optimization variables $\{x_n\}_{n=1}^N$ are also partitioned among nodes, and each variable is of interest to exactly one node. To this end, let $\{\S_k\}_{k=1}^K$ denote disjoint subsets such that the variables $\{x_n | x_n \in S_k\}$ are local to node $k$. Further, the component function $g_k(\cdot)$ at node $k$ depends only on variables that are local either to node $k$ or its neighbors. The overall problem considered here takes the following form:
\begin{align}\label{prob2}   
 P=\min_{\x}   & \sum_{ k=1 }^{ K } g_{ k }(\{x_n\}_{n \in \S'_k}) + h_k(\{x_n\}_{n\in \S_k}) \\  
 \text{ s. t. }  & \{x_n\}_{n\in \S_k} \in  \cX_k \quad k = 1, 2, \ldots, K. \nonumber
\end{align}
where the set $\S_k':=\bigcup_{j \in \{k \} \cup \N'_k} \S_j$. Observe here that for each node $k$, the function $h_k$ and the constraint set $\cX_k$ depend only on the variables in $\S_k$. { For this partially separable form, it is now possible to express the objective function as a bipartite factor graph, with check nodes representing the summands $g_k$, and the variable nodes representing the sets $\mathcal{S}_k$. For the four-node example network shown in Fig. \ref{network}, the factor graph is shown in Fig. \ref{separable:fig}. From the perspective of algorithm design, the dependence structure imposed by \eqref{prob2} can be exploited to eliminate message passing between non-neighboring nodes.}  The partially separable form considered here occurs commonly in the context of distributed estimation, where the parameter of interest is a collection of node-specific quantities such as temperatures, node locations,  harmful algal blooms, pH, and temperature \cite{cressie2015statistics,kekatos2013distributed}. Cooperation between nodes is still required here, since parameters at neighboring nodes are often coupled or correlated. The subsequent examples detail two specific applications where the { partially separable form} problem structure  in \eqref{prob2} arises. 

\subsection{Examples}\label{eg}

\indent \textbf{Example 1. }
Consider a sensor network with nodes $\cK:=\{1,2, \ldots, K\}$ observing an environmental phenomenon. The environmental state at the location of sensor $k \in \cK$ is denoted by the random variable $\theta_k \in \Rn$, and is collected in the vector $\bt \in \Rn^{K}$. The measurements at sensor $k$ are in form of a vector $\y_k \in \Rn^{N_k}$ and depend on the local state $\theta_k$ via the pdf $p(\y_k \mid \theta_k)$.  Given the overall measurement vector $\y$, that collects the nodal measurements $\{\y_k\}_{k=1}^K$, the goal at each node is to calculate the \emph{maximum a posteriori} (MAP) estimate of $\theta_k$, given by $\hat{\bt} := \arg \max_{\bt} p(\bt \mid \y) = \arg \max_{\bt} \ln p(\bt) + \sum_{k=1}^K \ln p(\y_k \mid \bt)$, where $p(\bt)$ captures the available prior information about $\bt$.  Of particular interest here are problems where the prior pdf has a factored representation, i.e., $p(\bt) = \prod_{k=1}^N p_k(\bt_{C_k})$, where $C_k \subset \cK$, and $\bt_{C_k}$ stacks the variables $\{\theta_n\}_{n \in C_k}$. Such factored pdfs arise naturally within the context of distributed sensing of Markov random fields \cite{dogandvzic2006distributed}, where $C_k  = \N_k$. The Markov property implies that $\theta_k$ is correlated with $\{\theta_n\}_{n \in \N_k}$, but is conditionally independent of $\{\theta_n\}_{n \notin \N_k}$, given $\{\theta_n\}_{n \in \N_k}$; see e.g\cite{rajawat2012network} and references therein. When the pdf can be factored, the MAP estimation problem may be written as
\begin{align}\label{map2}
\hat{\bt} = \arg \max_{\bt} \sum_{k=1}^K \ln p_{N_k}(\bt_{\N_k}) + \sum_{k=1}^K \ln p(\y_k \mid \theta_k).
\end{align}
It can be observed that \eqref{map2} is a special case of \eqref{prob2} with $g_k(\cdot) = \ln p_k(\cdot)$ and $h_k(\cdot) = \ln p(\y_k | \cdot)$.  

\indent \textbf{Example 2. } This example details the problem of cooperative localization in wireless networks via multidimensional scaling. Distributed localization is necessary in networks where manual or GPS-assisted positioning is not feasible, and not all nodes can communicate with the fusion center or base stations. Within the cooperative localization framework considered here, the nodes utilize pairwise range measurements and the known location of a few anchor nodes to estimate their locations. Range measurements are often made using techniques such as received signal strength or time-of-arrival \cite{rss_crlb}. Specifically, consider a network with $N$ nodes, distributed in a $p$-dimensional space with $p = 2$ or $3$. The location of node $k$ is denoted by $\x_k \in \Rn^p$, and collected into the matrix $\X=[ \x_1,\x_2,\ldots,\x_N ]$. The measured distance between a pair of neighboring nodes $k$ and $j$ is denoted by $\delta_{kj}$, and is generally a noisy version version of $d_{kj}(\x_k,\x_j) := \norm{\x_k-\x_j}$. The MDS formulation finds the overall configuration $\X$ by solving the following sum of squares problem \cite{localization_example, hero}
\begin{align}\label{mds}
\hat{\X}=\arg \min _{\X \in \mathbb{R}^{n\times p}} \sum_{k=1}^{N} g_k(\{\x_j\}_{j \in \N_k})
\end{align}
where, $g_k(\{\x_j\}_{j \in \N_k})=\sum_{j \in \N'_k } w_{kj}\left(\delta_{kj}-d_{kj}(\x_k,\x_j)\right)^2$. When the location $\x_k^p$ of node $k$ is known a priori, a regularization term of the form $r_i\norm{\x_k-\x^p_k}$ may also be added to $g_k(\cdot)$. Since \eqref{mds} is non-convex, it is often solved via majorization \cite{hero} or via SDP relaxation \cite{localization_example}. A distributed and incremental algorithm for solving \eqref{mds} via majorization was first detailed in \cite{hero}. A distributed and asynchronous algorithm using SDP relaxation was described in \cite{localization_example}. The subsequent sections describe the proximal and majorized ADMM algorithms that are also applicable to \eqref{mds}.


\section{Distributed Asynchronous ADMM}\label{asec1}
This section details the proposed distributed asynchronous algorithm for solving \eqref{prob2} via ADMM, and provides the relevant convergence results. To begin with, the next subsection describes a distributed synchronous implementation, which is motivated from the so called general-form consensus algorithm for solving convex optimization problems \cite{admm_defn}. The synchronous version serves as a starting point for the  asynchronous algorithm described in Sec. \ref{asec}. 

\subsection{Distributed Synchronous ADMM} \label{dsadmm}
For ease of exposition, it is assumed henceforth that the sets $\S_k$ are singletons, i.e., $\S_k = \{k\}$. Note however that the algorithms and the corresponding convergence results are applicable to the general case as well. For instance, it is possible to have $\abs{\S_k} = n_k$ with $\sum_k n_k = N$ and $g_k: \Rn^{N_k} \rightarrow \Rn$ where $N_k = n_k + \sum_{j \in \N'_k} n_j$. Indeed, for the localization example described in Sec. \ref{eg}, $n_k = 2$ for all $k$. The expanded neighborhood of node $k$ is defined as $\N_k := \{k\} \cup \N'_k$ and includes $k$ itself. Towards developing a distributed algorithm for solving \eqref{prob2}, introduce copies $x_{kj}$ of the variable $x_k$ corresponding to all nodes $j \in \N_k$ for all $k \in \cK$. To make sure that the copies of the a variable agree across nodes, also introduce the so called consensus variable $z_k$ for each neighborhood $\N_k$. The introduction of these extra variables amounts to reformulating \eqref{prob2} as 
\begin{align}\label{prob_formulation}
 & \min_{\{\x_k\}, \z} \sum _{ k=1 }^{ K }{ g_{ k }(\{ x_{kj}\}_{j\in \N_k})  } + h_k(z_k) \\
& \text{s. t. } \quad   x_{kj}=z_j, \quad  j \in \N_k, \quad k=1,\ldots,K \label{cons1}\\
& z_k \in \cX_k,  \qquad \qquad k=1,\ldots,K \nonumber
\end{align}
where $\x_k$ collects the variables required at node $k$, i.e.,  $[\x_k]_j = x_{kj}$ for $j \in \N_k$, and zero otherwise. 

The idea of introducing consensus variables, in order to make the updates separable in the optimization variables, is well known \cite[Chap. 5]{parikh2013proximal}. It is now possible to apply the ADMM method by associating dual variables $y_{kj}$ for each constraint in \eqref{cons1} and writing the augmented Lagrangian as
\begin{align}\label{lagrangian}
& L(\{\x_k\},\z,\{\y_k\})  = \sum_{ k=1 }^{ K } \Bigg( { g_{ k }( \x_{ k }) } + h_k(z_k )\nonumber \\
 & + \qquad \sum _{ j \in \N_k }{ \langle y_{ kj }, x_{ kj }-z_j \rangle  } 
  +\sum _{ j \in \N_k  }{ \frac { \rho _{ k } }{ 2 } \norm{  x_{ kj }-z_j } ^{ 2 } } \Bigg)
\end{align}
where $\rho_k >0$ is a positive penalty parameter and $ \z$ collects $\{z_k\}_{n \in \cK}$. Here, although $\{x_{kj}\}$, $\{z_j\}$, and $\{y_{kj}\}$ are all scalar variables, the notations for generalized inner product and norm are retained in order to avoid manipulations that are specific to the scalar case. On the left-hand side, for each $k$, the dual vector $\y_k \in \Rn^N$ is such that $[\y_k]_j = y_{kj}$ if $j \in \N_k$, and zero otherwise.

From \eqref{lagrangian}, it is clear that $L$ is separable in $\{\x_k\}$. The Lagrangian is also separable in $\{z_j\}$, since it holds that

\begin{align}\label{lang2}
 L(\{ \x_{ k }\},\z, \{\y_k\}) 	& =\sum_{ j=1 }^{ K }\Bigg(g_{ j }(  \x_{ j })  + \sum _{ k \in \N_j }{ \langle y_{ kj }, x_{ kj }-z_j \rangle  } \nonumber \\ & +\sum _{ k \in \N_j  }{ { \frac { \rho _{ k } }{ 2 } (  x_{ kj }-z_j ) ^{ 2 } } }+h_j(z_j)  \Bigg). 
\end{align}

Together, \eqref{lagrangian} and \eqref{lang2} allow the ADMM updates to be carried out in a distributed fashion. In particular, starting with arbitrary $\{\x_k^1\}$ and $\{\y^1_{kj}\}$, the update for $\{z_j^{t+1}\}$ are evaluated as
\begin{align}
z_j^{t+1} & = \arg \min _{ z_j\in \cX_j } h_j(z_j) + \sum _{ k \in \N_j }{ \langle y_{ kj }^t, x_{ kj }^t-z_j \rangle  } \nonumber \\ & \qquad +\sum _{ k \in \N_j  }{ { \frac { \rho _{ k } }{ 2 } \norm{  x^t_{ kj }-z_j } ^{ 2 } } }\nonumber\\
& = \pj\left( \frac{\sum_{k\in \N_j}\left(\rho_k x_{kj}^t + y_{kj}^t\right)}{\sum_{k \in \N_j}\rho_k}  \right)\label{z_k_update}
\end{align}
where the proximal point function $\pj(\cdot)$ is defined as
\begin{align}
\pj(x) := \arg \min_{u \in \cX_j} h(u) + \frac{1}{2}\norm{x-u}^2.
\end{align}

Similarly updates for $\x_k^{t+1}$ can be obtained by minimizing \eqref{lagrangian} with respect to $\x_k$, which yields 
\begin{align}
\x_{k}^{t+1}= & \arg\min _{ \x_k} { g_{ k }( \x_{ k }) }  +\sum _{ j \in \N_k }{ \langle y_{ kj }^{t}, x_{ kj }-z_j^{t+1} \rangle  } \nonumber \\ &  \quad + \sum _{ j \in \N_k  }{ \frac { \rho _{ k } }{ 2 } \norm{  x_{ kj }-z_j^{t+1} } ^{ 2 } } \label{x_k_sync}
\end{align}
where the optimization is with respect to $\{x_{kj}\}_{j \in \N_k}$. Observe that since the component functions $g_k(\cdot)$ are not necessarily convex, the update in \eqref{x_k_sync} is difficult to carry out. As suggested in \cite{admm}, the update $\x_k^{t+1}$ can however be calculated approximately as follows
\begin{align}\label{x_k_sync1}
& \x_{k}^{t+1} \approx \arg \min _{ \x_k } { g_{ k }( {\z}_{ k }^{t+1}) }+ \langle \nabla g_{ k }( {\z}_{ k }^{t+1}), \x_{ k }-{\z}_{ k }^{t+1} \rangle \nonumber \\ &
 + \sum _{ j \in \N_k }{ \langle y_{ kj }^t, x_{ kj }-z_j^{t+1} \rangle  }+\sum _{ j\in \N_k  }{ { \frac { \rho _{ k } }{ 2 } \norm{  x_{ kj }-z_j^{t+1} } ^{ 2 } } }
\end{align}
where the vector $[{\z}_k]_j := z_j$ for all $j \in \N_k$ and zero otherwise. Since nodal functions $g_k(\cdot)$ depend only on $\{x_n\}_{n \in \N_k}$, the gradient vector is defined as
\begin{align}\label{grad_syn}
[\nabla g_{ k }( {\z}_{ k }^{t+1})]_j := \begin{cases} \frac{\partial}{\partial x_{kj} }g_k(\x_k) \Big|_{\x_k=\z_k} &  j \in \N_k \\
0 & j \notin \N_k.\end{cases}
\end{align}
The approximate update of $x_{kj}^{t+1}$ thus becomes
\begin{align}\label{x_k_syn}
x_{ kj }^{ t+1 }= \begin{cases} z_j^{t+1} - \frac{1}{\rho_k}\left([\nabla g_{ k }( {\z}_{ k }^{t+1})]_j+y_{kj}^{t}\right) &  j\in \N_{ k }  \\
 0 & j \notin \N_k. \end{cases}
\end{align}
Finally, dual updates are given by 
\begin{align}\label{eq_y_k1}
 y_{ kj }^{ t+1 } &= y^{ t }_{ kj }+\rho _{ k }\{ x_{ kj }^{ t+1 }-z_{ j }^{ t+1 }\} & j \in \N_k
\end{align}

Algorithm \ref{sync_algo1} summarizes the implementation of the distributed ADMM described here. The main feature of Algorithm \ref{sync_algo1} is that it does not require a master node, and all message passing is limited to the neighboring nodes only. The distributed implementation also requires that each node must be capable of carrying out the updates in \eqref{z_k_update} and \eqref{x_k_sync1}. Algorithm \ref{sync_algo1} may be viewed as the proximal variant of the distributed ADMM algorithm \cite{admm_defn}, applied to  the non-convex problem \eqref{prob2}. The stopping criterion in Algorithm \ref{sync_algo1} is simply $\abs{z_k^{t+1}-z_k^t} \leq \delta$ for all $k \in \cK$ and a small $\delta>0$.

\begin{algorithm} 
	\caption{: Distributed  Algorithm}
	\begin{algorithmic}[1]\label{sync_algo1}
		\STATE Initialize $\{ x_{kj}^1,y_{kj}^1 \}$, $\z_k$ for all $j \in\N_k$.
		\FOR{$t=1,2,\ldots$}
		\STATE Send $\{\rho_k x_{kj}^{t}+y_{kj}^{t} \}$ to neighbors $j \in \mathcal{N}_k$
		\STATE  Upon receiving $\{\rho_j x_{jk}^{t}+y_{jk}^{t} \}$ for all $j\in\N_k$,
		\STATE \hspace{1cm} Update $z_{k}^{t+1}$ as in \eqref{z_k_update}
		\STATE Transmit $z_k^{t+1}$ to its neighbors $j\in \N_k$
		\STATE Upon receiving $z^{t+1}_j$, from all neighbors $j\in\N_k$, 
		\STATE \hspace{1cm} Update the primal variable $\x_k^{t+1}$ as in \eqref{x_k_syn}
		\STATE \hspace{1cm} Update the dual variable $y_{kj}^{t+1}$ as in \eqref{eq_y_k1}
		\IF{stopping criteria is met}\STATE terminate loop \ENDIF
		\ENDFOR
	\end{algorithmic}
\end{algorithm}

Algorithm \ref{sync_algo1} is a synchronous protocol since all updates must necessarily be carried out at every iteration by every node. Its applicability to heterogeneous networks is therefore limited, since the progress of the algorithm is determined by the slowest node in the network. For instance, the following issues may arise when Algorithm \ref{sync_algo1} is implemented on a wireless network with energy-constrained, low-cost devices. 
\begin{enumerate}
	\item[\textbf{(S1)} ] For some nodes, calculating $\nabla g_{k}(\cdot)$ [cf. step 8] or $\pk(\cdot)$ [cf. step 5] may be computationally demanding. In such cases, all nodes in the network must wait for the slowest node to carry out its update. 
	\item[\textbf{(S2)} ] Each node is required to transmit two messages to each of its neighbors per-iteration. This might be excessive for nodes operating on a power-budget.

\end{enumerate}

Delays in function computation or communication may also arise because of the heterogeneity of nodes in a wireless network. For instance, mission critical nodes may attempt to extend their battery lives by operating under a power-saving mode, thus deliberately reducing their computational capabilities. On the other hand, the available power at some energy-harvesting nodes may vary throughput the day, depending on, for instance, the received solar energy. 

Sec. \ref{asec} describes an asynchronous version of Algorithm \ref{sync_algo1} that overcomes \textbf{(S1)}-\textbf{(S2)} by allowing nodes to carry out the resource-intensive operations in Steps 5 and/or 7 intermittently. Each time an update is skipped, the node saves on both computational and communication costs. It is shown that the asynchronous algorithm converges, as long as the updates are performed ``often enough,'' a notion that is made precise through some constraints on the algorithm parameters.

\subsection{Distributed Asynchronous Algorithm with Optional and Delayed Updates}\label{asec}
This section describes the asynchronous distributed ADMM algorithm that allows nodes to skip the calculation of $\nabla g_{k}(\cdot)$ [cf. step 8] or $\pk(\cdot)$ [cf. step 5] for some iterations. For step 5, this is achieved by simply setting $z_{j}^{t+1} = z_j^{t}$ for each non-updating node $j \in \cK$ at time $t$. Define $\S^t$ as the set of nodes that carry out the update at time $t$, then Step 5 becomes
\begin{align}\label{zasynch}
z_j^{t+1} &= \begin{cases} \pj\left( \frac{\sum_{k\in \N_j}\left(\rho_k x_{kj}^t + y_{kj}^t\right)}{\sum_{k \in \N_j}\rho_k}  \right) & j \in \S^t \\
z_j^{t} & j \notin \S^t
\end{cases}
\end{align} 
Whenever $j \notin \S^t$, the subsequent transmission may not be carried out either, and the non-updating node may simply stay silent. The neighboring nodes will then wait for a fixed amount of time to receive an update, and assume that $z_j^{t+1} = z_j^t$ holds for all nodes $j \in \N_k$ that do not transmit anything. 

For step 7, each node uses the latest available gradient $\nabla g_k (\z_k^{{[t+1]}_k})$ for update in \eqref{x_k_syn}, where $t+1-T_k \leq {[t+1]}_k \leq t+1$ for some $T_k < \infty$. In other words, the gradient calculation may only be carried out intermittently, and the same gradient can be used for next several time slots. Alternatively, for computationally challenged nodes, the gradient calculation itself may take several time slots. Dropping the subscript $k$ from $[t+1]$ for notational brevity, the update in Step 7 becomes
\begin{align}\label{xasynch}
x_{ kj }^{ t+1 }= \begin{cases} z_j^{t+1} - \frac{1}{\rho_k}\left([\nabla g_{ k }( \z_{ k }^{[t+1]})]_j+y_{kj}^{t}\right) &  j\in \N_{ k }  \\
 0 & j \notin \N_k. \end{cases}
\end{align} 
The transmission of updates in Step 3 is again optional and in the event that no update is received at a neighbor $j \in \N_k$,  \eqref{zasynch} is used instead.

\begin{algorithm}

	\caption{Distributed  Asynchronous ADMM with Optional Updates}
	\begin{algorithmic}[1]\label{async_algo1}
		\STATE Set $t=1$, initialize $\{ x_{kj}^1,y_{kj}^1,z_j^1\}$ for all $j \in\N_k$.
		\FOR{$t=1,2,\ldots$}
		\STATE (Optional) Send $\{\rho_k x_{kj}^{t}+y_{kj}^{t} \}$ to neighbors $j \in \mathcal{N}_k$ \label{sendxa1}
		\IF{$\{\rho_j x_{jk}^{t}+y_{jk}^{t} \}$ received from all $j\in\N_k$} \STATE{(Optional) Update $z_{k}^{t+1}$ as in \eqref{z_k_update} and transmit to each $j \in \N_k$} \label{sendza1} \ENDIF
		 \IF{$z^{t+1}_j$ not received from some $j \in \N_k$} \STATE{set $z_j^{t+1} = z_j^t$} \ENDIF 
		\STATE (Optional) Calculate gradient $\nabla g_{ k }(\z_{ k }^{t+1})$
		\STATE Update the primal variable $\x_k^{t+1}$ as in \eqref{xasynch} \label{xupa1}
		\STATE Update the dual variable $y_{kj}^{t+1}$ as in \eqref{eq_y_k1} \label{yupa1}
		\IF{$\norm{\x_k^{t+1}-\x_k^t} \leq \delta$ \label{stopca1}}\STATE terminate loop \ENDIF
		\ENDFOR
	\end{algorithmic}
\end{algorithm}

The proposed asynchronous algorithm for node $k \in \cK$ is summarized in Algorithm \ref{async_algo1}.  The salient features of the proposed asynchronous algorithm are as follows.
\begin{enumerate}
	\item All resource-intensive steps, such as gradient calculation, proximal function calculation, and transmission of updates are now optional.
	\item Nodes operating on a power budget may only carry out the updates in Steps \ref{xupa1} and \ref{yupa1} at every iteration. When using the old gradient, these updates amount to simple addition/subtraction operations. It is also possible to defer Steps \ref{xupa1} and \ref{yupa1} to a later time slot when the transmission in Step \ref{sendxa1} occurs.
	\item The nodes may implement a timeout mechanism when listening for updates [cf. Steps \ref{sendxa1} and \ref{sendza1}]. The appropriate default actions must be triggered if nothing is heard from one or more neighbors. 
\end{enumerate}

\begin{figure}[tbh]
	\centering
	\includegraphics[scale=.50]{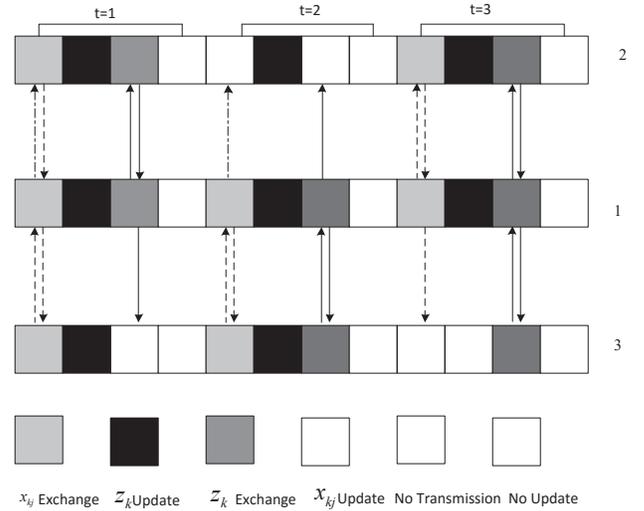}
	\caption{Timing Diagram}
	\label{timing}
\end{figure}

Fig. \ref{timing} shows an illustration where nodes 2 and 3 are neighbors of a node 1, that go to sleep at iteration $t=2$. During the sleep state, a node is only able to receive updates, but not transmit or carry out any computations. Observe that in the second sub-slot of the 2-nd iteration, $z^3_2$ and $z^3_3$ are also not updated at nodes $2$ and $3$ since no update for $x_{12}^{2}$ and $y_{12}^2$ was received. Consequently, it will hold that $z_k^{t} = z_k^{t-1}$ for  $k = 1$, $2$, and $3$. More generally, each sleeping node will force all its neighbors to use the update in Step 8 of Algorithm \ref{async_algo1}. Since the update frequency at each node cannot be too small, it is required that several nodes should be awake at any given time.

The next section provides the convergence analysis of Algorithm \ref{async_algo1}, ensuring that for any $\delta > 0$, the stopping criteria in Step \ref{stopca1} of Algorithm \ref{async_algo1} is eventually met. Note that the convergence results also apply to Algorithm \ref{sync_algo1}, which is simply a special case of Algorithm \ref{async_algo1}. Before concluding this subsection, the following remark about related work in the context of asynchronous algorithms is due. 

\begin{rem}
The key feature of the proposed algorithms is that some of the updates and transmissions are entirely optional. From the perspective of algorithm design, this results in two sources of asynchrony, namely, delayed gradients, and skipped updates. The asynchronous algorithm in \cite{admm} can handle bounded delays in the gradient calculation, but does not consider skipping updates. This is because the DA-ADMM algorithm in \cite{admm} utilizes a fusion center that is not energy constrained. Note that the absence of such a fusion center significantly complicates the algorithm design, also restricting the application of Algorithms \ref{async_algo1} to \eqref{prob2}, which is a special case of \eqref{prob1}. 

An asynchronous ADMM algorithm was also proposed and applied to non-convex problem in \cite{async_admm}. However, no proof of convergence was provided, and the performance of the algorithm was only tested via simulations. A number of asynchronous variants exist for convex problems. The optional-update idea used in \eqref{zasynch} is in fact inspired from the asynchronous ADMM proposed in \cite{chang2014proximal,order}. Different from Algorithm \ref{async_algo1} however, the asynchronous algorithm in \cite{chang2014proximal} also uses random network states, and is applicable only to convex problems. For partially separable problems such as \eqref{prob2}, it may also be possible to develop an asynchronous variant of the coordinate descent algorithm. The stochastic asynchronous coordinate descent algorithm proposed in \cite{liu2013asynchronous} also allows erroneous gradients but is applicable only to convex problems. Finally, asymptotic results on the effect of asynchrony on the stochastic gradient descent algorithm were provided in \cite{lian2015asynchronous}. 
\end{rem}
  
\subsection{Convergence Analysis for Algorithm \ref{async_algo1}}\label{conv1}
In order to establish the convergence of the asynchronous algorithm, some assumptions regarding the problem structure and algorithm parameters must be made. In particular, it is shown that the algorithm converges as long as the optional updates happen ``often enough.'' Specifically, recall that $t+1-[t+1] \leq T_k$, which implies that the gradients may be calculated using updates that are at most $T_k$-old. Similarly, let the frequency of update in \eqref{z_k_update} being carried out at node $k$ be denoted by $0 < f_k \leq 1$. For instance, if the update occurs once every $K$ time slots, $f_k = 1/K$. Then the following assumptions are required.

\begin{assumption}  \label{aslip}	
For each node $k$, the component function gradient $\nabla g_k(\x)$ is Lipschitz continuous, that is, there exists $L_k > 0$, for all $\x,\x' \in \mathrm{dom} g_k$ such that
  		\begin{align}
  		\left\| \nabla g_k( \x) - \nabla g_k(\x') \right\| \le L_k\left\| \x-\x' \right\|.
  		\end{align}
		  
\end{assumption}
\begin{assumption} \label{asX} The set $ \mathcal{X}$ is a closed, convex, and compact. The functions $g_k(x)$ is bounded from below over $\mathcal{X}$.
\end{assumption}

\begin{assumption} \label{asrho}
For node $k$, the step size $\rho_k$ is chosen large enough such that, it holds that $  \alpha_k>0$ and $\beta_k>0$, where
  	 \begin{align}
  \alpha_k &:= \frac{\rho_kf_k}{2}-  \left(   \frac { 7 
  	 	L_k }{ 2 \rho _{ k }^2 }+\frac{1}{\rho_k}\right)|\N_k|L_k^2(T_k+1)^2-\frac{|\N_k|L_kT_k^2}{2}    \nonumber\\
  	\beta_k &:= \rho_k - 7 L_k \label{rholk}
  	 \end{align}
  \end{assumption}

Of these, Assumptions (A1) and (A2) are standard in the context of non-convex optimization \cite{hong2014convergence,admm,magnusson2014convergence} and are satisfied for most problems of interest. Intuitively, the use of old gradients is permitted only because they are Lipschitz continuous, and therefore change slowly over iterations. The boundedness assumption (A2) is required to ensure that the updates $\{\x^{t},\z^t,\y^t\}$ generated by Algorithm \ref{async_algo1} stay bounded. This assumption may be dropped for certain problems where boundedness of these iterates may arise naturally. Finally, Assumption (A3) specifies the exact relationship that the algorithm and problem parameters $\{L_k,\rho_k,f_k,T_k\}_{k=1}^K$ must satisfy in the worst case. Interestingly, the choice of $\rho_k$ is node-specific since nodes may differ with respect to their computational capabilities resulting in different levels of asynchrony. While not stated explicitly, the convergence also requires that $\rho_k < \infty$. Thus, for \eqref{rholk} to hold, it is necessary that $f_k > 0$ and $T_k < \infty$. In other words, the delay cannot be unbounded in the worst case. 

The convergence of the asynchronous algorithm is established through the following intermediate lemma that holds under Assumptions (A1)-(A3).

\begin{lem}\label{lem1}
\begin{enumerate}
\item[(a)] Starting from any time $t=t_0$, there exists $T < \infty$ such that
\begin{align} \label{lagdiff}
 & L( \{\x_{ k }^{ T+t_0 }\}; {\z}^{T+t_0},\{\y_k^{ T+t_0 }\})-L( \{\x_{ k }^{ t_0 }\};   {\z}^{t_0},\{\y_k^{ t_0 }\})\nonumber\\
 &\le -\sum_{i=t_0}^{T+t_0-1}\sum_{k=1}^{K}\frac{\beta_k}{2}\sum_{j\in \N_k}\|x_{kj}^{i+1}-x_{kj}^{i}\|^2 \nonumber \\
 & \hspace{1cm}-\sum_{i=t_0}^{T+t_0}\sum_{k=1}^{K}\alpha_k \sum_{j\in \N_k}\|z_{j}^{i+1}-z_{j}^{i}\|^2.
 \end{align}
 \item[(b)] The augmented Lagrangian  {   values in \eqref{lang2}  }  are bounded from below, i.e., for any time $t \geq 1$, it holds that
 Lagrangian satisfies
 \begin{align}
 L(\{ \x_k^{t}\}; {\z}^{t},\{\y_k^{ t }\}) \geq \mathsf{P} - \frac{L_k}{2}\sum_{j \in \N_k} \mathrm{diam}^2(\cX) >  -\infty \nonumber
 \end{align}
\end{enumerate}
\end{lem}

The proof of Lemma \ref{lem1} is provided in Appendix A. Lemma \ref{lem1}(a) establishes that there exists some finite $T$ such that the augmented Lagrangian values are non-increasing after $T$ iterations. In practice, the value of $T$ depends on the update frequencies $\{f_k\}_{k=1}^K$. For instance, $T$ could be the minimum number of iterations in which each node $k$ updates $z_k^t$ at least $f_kT$ times. Lemma \ref{lem1}(b) establishes that the Lagrangian is bounded from below. One way to interpret Lemma \ref{lem1} is to define the sequence $L^{n}:= L( \{\x_{ k }^{ nT+1 }\}; {\z}^{nT+1},\{\y_k^{ nT+1}\})$ for all $n \geq 0$, and observe that $\{L^n\}$ is non-increasing and bounded from below, and therefore convergent. The subsequent theorem establishes the final convergence result and related properties.

\begin{theorem}\label{thm1}
\begin{enumerate}
 \item[(a)] The iterates generated by Algorithm \ref{async_algo1} converges in the following sense
\begin{subequations}
\begin{align}\label{th1}
\lim_{ t\rightarrow \infty  }{\norm{z_{k}^{t+1} -   z_k^{t} }  } & =0, \quad \forall \quad k \\
\lim_{ t\rightarrow \infty  }{\norm{x_{kj}^{t+1} -   x_{kj}^{t} }  }  &=0,  \quad j  \in \N_k, \quad\forall \quad  k\\
\lim_{ t\rightarrow \infty  }{\norm{y_{kj}^{t+1} -   y_{kj}^{t} }  } & =0, \quad j  \in \N_k, \quad \forall \quad k
\end{align}
\end{subequations} 
\item[(b)] For each $k \in \cK$ and $j \in \N_k$, denote limit points of the sequences $\{z_k^{t}\}$, $\{x_{kj}^{t}\}$, and $\{y_{kj}^{t}\}$ by $z_k^{\star}$, $x_{kj}^{\star}$, and $y_{kj}^{\star}$, respectively. Then $\{\{z_k^\star\}, \{x_{kj}^\star\}, \{y_{kj}^\star\}\}$ is a stationary point of \eqref{prob_formulation} and satisfies
\begin{subequations}
\begin{align}
&\nabla g_k(\x_k^\star)+ \y_k^\star = 0,\quad  k=1,\ldots,K \label{th2_a}  \\ 
& \sum_{j \in \N_k}y_{jk}^\star \in \partial(h_k(z))\mid_{z = z_k^\star} \quad k=1,\ldots,K \label{th2_b}\\ 
& \quad x_{kj}^\star=z_j^\star \in \cX_j,  \quad  j\in \N_k,  k=1,\ldots,K \label{th2_c} 
\end{align}
\end{subequations}
\end{enumerate}
\end{theorem}

The proof of Theorem \ref{thm1} is provided in Appendix \ref{{theorem}_proof}. Note that it suffices to show that Algorithm \ref{async_algo1} converges to a stationary solution of \eqref{prob_formulation} which is equivalent to  \eqref{prob2}. In other words, $\{z_n^\star\}_{n=1}^N$ can be used as a solution to \eqref{prob2}. It is emphasized that Algorithm \ref{async_algo1} may not necessarily converge to a globally optimum solution to \eqref{prob2}.

 Further, using assumptions (A1)-(A3) alone, it is also difficult to quantify the convergence rate of the ADMM algorithm for the non-convex case. As will be shown in Sec. \ref{sim} however, rate of convergence  for the localization problem is low whenever $\rho_k$ is large. Assuming that this observation applies to \eqref{prob2} generally, it makes sense to choose $\rho_k$ as small as possible, while respecting  \eqref{rholk}. Interestingly, this result also matches with the intuition that the Algorithm \ref{async_algo1} converges slowly if the asynchrony is high, i.e., when $T_k \gg 1$ and $f_k \ll 1$, since it would require choosing a larger $\rho_k$ for each $k \in \cK$. 

\section{Majorized Asynchronous ADMM}\label{madmms}
In many problems, it is possible to upper bound the non-convex component functions $g_k(\x_k)$ with an appropriate convex surrogate function. Given $\z_k \in \cX$, the surrogate $f_k(\x_k,\z_k)$, also referred to as the majorizing function, is such that for all $\x_k$, $g_k(\x_k) \leq f_k(\x_k,\z_k)$, and satisfies
\begin{align}
f_k(\z_k,\z_k) &= g_k(\z_k) \label{geq} \\
\nabla f_k(\z_k,\z_k) &:= \nabla_{\x} f_k(\x,\z_k) \mid_{\x = \z_k} = \nabla_{\x} g_k(\x)  \mid_{\x = \z_k} .\label{gradeq}
\end{align}
When such a majorizer exists and can be found easily, the update of $\x_k^{t+1}$ in \eqref{x_k_syn} can be carried out more accurately, without appreciable increase in the per-iteration complexity. 

This section develops a provably convergent variant of Algorithm \ref{async_algo1} that utilizes majorization while updating $\x_k^{t+1}$. Specifically, while \eqref{eq_y_k1} and \eqref{zasynch} stay the same, the update for $\x_k^{t+1}$ becomes
\begin{align}\label{mxasynch}
\x_k^{t+1} =\arg\min f_k(\x_k, &\z_k^{[t+1]}) + \sum_{j \in \N_k} \langle y_{kj}^t, x_{kj}-z^{t+1}_j \rangle + \nonumber \\
& \frac{\rho_k}{2}\sum_{j \in \N_k} \norm{x_{kj}-z_j^{t+1}}^2
\end{align}
where, as in \eqref{xasynch},  $\z_k^{[t+1]}$ is used for calculating the surrogate function. Interestingly, the majorized ADMM proposed here enjoys the flexibility afforded by its proximal counterpart, and can be implemented as in Algorithms \ref{async_algo1} or \ref{async_algowsn}.

In order to show that the majorized ADMM converges, the following assumptions are required in addition to (A2).
\begin{assumption}  \label{mlip}	
For each node $k$, there exists a constant $L_k \geq 0$, such that for all $\x,\x',\z,\z'$, the following inequalities are satisfied:
\begin{subequations}\label{fglip}
\begin{align}
\norm{ \nabla g_k( \x) - \nabla g_k(\x')}  &\leq L_k\norm{ \x-\x'} \label{l1}\\
\norm{\nabla f(\x,\z) - \nabla f(\x',\z)} & \leq L_k\norm{\x-\x'} \label{l2}\\
\norm{\nabla f(\x,\z) - \nabla f(\x,\z')} & \leq L_k\norm{\z-\z'}. \label{l3} 
\end{align}
\end{subequations}
\end{assumption}

\begin{assumption} \label{asrho2}
For node $k$, the step size $\rho_k$ is chosen large enough such that $\alpha_k>0$ and $\beta_k>0$, where
	\begin{align}
	\alpha_k &:=  |\N_k|\left[\left(\frac{\rho_kf_k}{2}-  \left(   \frac { 8 
		L_k }{\rho _{ k }^2  
	}+\frac{1}{\rho_k}\right)L_k^2(T_k+1)^2-\frac{L_kT_k^2}{2}\right) \right] \nonumber  \\
	\beta_k &:= \frac { \rho _{ k }-9L_k}{ 2 }-\frac{8L_k^3}{\rho_k^2}. \label{betak2}
	\end{align}
  \end{assumption}

Note that Assumption \eqref{mlip} utilizes the same Lipschitz constant in \eqref{fglip} for the sake of simplicity. In general, if there exist constants $L^g_k$, $L^x_k$, and $L^z_k$ in \eqref{fglip}, it is always possible to define $L_k:=\max\{L^g_k,L^x_k,L^z_k\}$ for all $k$. Likewise, the convergence proof provided here can also be developed for the case when the three constants are different, resulting in slightly tighter bounds.

The conditions required in \eqref{mlip} are not very restrictive. Consider for instance a component function $g_k(\x)$ that is expressible as a sum of a convex function $g^1_k(\x)$ and  a concave function $g^2_k(\x)$. The concavity of $g^2_k(\x)$ allows it to be majorized by its supporting hyperplane, i.e., given any $\z \in \textrm{dom}~g^2_k(\x)$,
\begin{align}
g_k(\x) & = g^1_k(\x) + g_k^2(\x) \nonumber \\ & \leq f_k(\x,\z)
:= g^1_k(\x) +  g^2_k(\z) + \langle \nabla g_k^2 (\z), \x-\z\rangle
\end{align} 
where it can be verified that $f_k(\x,\z)$ satisfies \eqref{geq} and \eqref{gradeq}. Then, observe that $\nabla f_k(\x,\z) = \nabla g_k^1(\x) + \nabla g_k^2(\x)$ for all $\x$ and $\z$. If $L^1_k$ and $L^2_k$ are Lipschitz constants of $\nabla g^1_k$ and $\nabla g^2_k$ respectively, it may be seen that \eqref{l1} and \eqref{l2} hold with Lipschitz constant $L^1_k+L^2_k$, while the left-hand side of  \eqref{l3} is identically zero. The following Theorem, whose proof is provided in Appendix \ref{madmm}, summarizes the main result of this section.

\begin{theorem}\label{thm2}
The iterates generated by the majorized ADMM (\eqref{zasynch}, \eqref{mxasynch}, and \eqref{eq_y_k1}) converge to a stationary point of \eqref{prob_formulation}. 
\end{theorem}

\section{Distributed Localization in Networks}\label{sim}
This section builds upon the localization examples introduced in Sec. \ref{eg} and provides simulation results for comparing two different related algorithms. Monte-Carlo simulations are performed over randomly generated networks with $N = 25$ nodes. To this end, nodes are uniformly distributed over a unit two-dimensional area $\mathcal{R}=[0,1]^2$. Of these, $m = 5$ nodes are located at $(0.25, 0.25)$, $(0.75, 0.25)$, $(0.25.0.75)$, $(0.5, 0.5)$, at $(0.75,0.75)$, and serve as anchor nodes for others. Distance measurement and communication between nodes $k$ and $j$ is possible only if they are within a distance of $R = 0.5$ of each other. For neighboring nodes, the weights $w_{kj}$ are all set to unity.

It is remarked that the classical SMACOF algorithm used for solving in \cite{mdsbook} is not distributed, and is therefore not applicable in the present context. The distributed weighted MDS (DwMDS) framework introduced in \cite{hero} is an incremental algorithm that utilizes component wise majorization in order to circumvent the non-convexity of the objective function. Different from the proposed algorithms, incremental algorithms utilize a cyclic message passing routine and are synchronous in nature. 

A distributed and asynchronous algorithm utilizing semidefinite relaxation (SDR) was proposed in \cite{localization_example}. The so-called E-ML (edge relaxed maximum likelihood) algorithm utilizes relaxation to first convert the objective into a convex one, and then applies distributed ADMM. Since ADMM is applied on a convex and relaxed version of \eqref{mds}, it can be made distributed and asynchronous using known techniques. Specifically, E-ML activates only a subset of edges $\cE^{(t)} \subset \cE$ at each iteration, allowing nodes to save energy by not transmitting on some links \cite{iutzeler2013asynchronous}. Different from the present Algorithms \ref{async_algo1}, which converges to a stationary point of \eqref{mds}, the E-ML algorithm converges to the global optimum of the relaxed version of \eqref{mds}. Since both limit points may be different from the global optimum of \eqref{mds} itself, the subsequent comparisons utilize the normalized root mean square error (NRMSE) of the location estimates, given by
$
NRMSE:=\sqrt{{\mathbb{E}\left[\|\hat{\X}-\X^\star\|_F^2\right]}/{\mathbb{E}\left[\|\X^\star\|^2_F\right]}}
$
where the expectations are performed via Monte-Carlo simulations using 100 independent realizations of the random network. 

In the present case, the following modified version of the objective function in \eqref{mds} is considered $ d_{kj}(\x_k,\x_j) = \sqrt{\norm{\x_k-\x_j} + \epsilon}$, where $\epsilon > 0$ is a small number introduced to make the objective function differentiable everywhere. This modification also makes $\nabla g_k$ Lipschitz continuous, as required in (A1), and verified in Appendix \ref{lips}. Recall that Algorithm \ref{async_algo1} allows two modes of asynchrony, namely, old gradient and skipped updates. In the implementation described here, it is assumed that at each time instant, some nodes go to sleep and do not carry out any update. Further, at each iteration, each node randomly chooses an older available gradient or calculates a more recent one, while ensuring that the gradient used is at most $T_k$-old.

\begin{figure}[tbh]
	\centering
	\includegraphics[scale=.6]{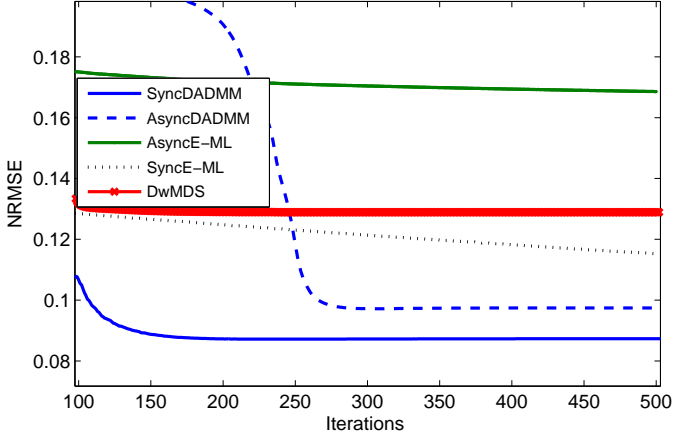}
	\caption{NRMSE performance}
	\label{localization}
\end{figure}

\subsection{NRMSE Performance}
Fig. \ref{localization} compares the NRMSE performance of DwMDS algorithm \cite{hero}, the synchronous (SyncE-ML) and asynchronous (AsyncE-ML) versions of the E-ML algorithm \cite{localization_example}, and the Algorithms \ref{sync_algo1} (SyncDADMM), \ref{async_algo1} (AsyncDADMM) . For the asynchronous E-ML algorithm, 125 of the 130 edges are activated per iteration. For Algorithm \ref{async_algo1}, we set $T_k = 8$ and $f_k = 0.75$, while $\rho_k$ is chosen so as to satisfy \eqref{asrho}. 

From Fig. \ref{localization}, it is clear that both, Algorithms \ref{sync_algo1}, \ref{async_algo1} outperform all the other state-of-the-art algorithms. As expected, the asynchronous versions of the E-ML and ADMM algorithms have poorer performance than their synchronous counterparts. Interestingly, asynchrony has a greater effect on the asymptotic performance of the E-ML algorithm than the convergence rate of the ADMM algorithm. It is also worthwhile to compare the implementation complexity of the two asynchronous algorithms. Since the E-ML algorithm uses SDR, its computational complexity is approximately $\mathcal{O}(\abs{\N_k}^3)$ for node $k$ per iteration. This is because each iteration requires solving an $\abs{\N_k}$-sized convex semidefinite program. On the other hand, gradient calculation is the most computationally demanding step in Algorithms \ref{sync_algo1} and \ref{async_algo1}, which translates to a total computational complexity of  $\mathcal{O}(\abs{\N_k})$ for node $k$ per iteration. Note further that unlike Algorithm \ref{async_algo1}, the updates in AsyncE-ML algorithm are also not optional, irrespective of the number of edges activated per iteration. The communication complexity for the two algorithms is of the same order, i.e., $\mathcal{O}(\abs{\N_k})$ for node $k$ per iteration. However, the message transmitted to each neighbor in AsyncE-ML consists of nine real numbers, as opposed to four real numbers in Algorithm \ref{async_algo1}.

\begin{figure}[tbh]
			\centering
			\includegraphics[scale=.5]{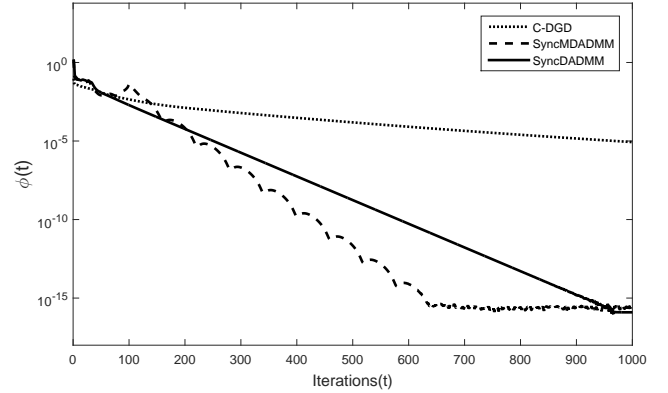}
			\caption{Convergence of Majorized and Proximal Variant and C-DGD}
			\label{majo prox cdgd}
\end{figure}

{

\subsection{Convergence rates}
As stated earlier, while Theorem \ref{thm1} establishes convergence, it does not provide any indication regarding the convergence rate exhibited by Algorithm \ref{async_algo1}. This subsection compares the convergence rates of the synchronous  proximal ADMM (SyncADMM), syncrhonous majorized ADMM  (SyncMDADMM), and the consensus-based distributed gradient descent (C-DGD) method proposed in \cite[Chap. 10]{palomar2010convex}. While C-DGD has only been proposed for convex problems, it has been implemented in a distributed and asynchronous manner \cite{nedic2011asynchronous, barbarossa2013distributed}, and is therefore an interesting candidate for the non-convex localization problem. The goal is to compare the three algorithms in terms of how fast they approach a local minimum, while ignoring their distances to the global minimum. To this end, Fig. \ref{majo prox cdgd}  shows the evolution of the convergence criterion $\phi(t):=\norm{\frac{1	}{N}\sum_{i=1}^N({\x_i^{t+1}-\x_i^t})}_F$ with the iteration index $t$. It can be observed that the majorized ADMM performs slightly better than the proximal ADMM. On the other hand, the C-DGD algorithm convergence relatively slowly. In particular, unlike the other simulations, the plots shown in Fig. \ref{majo prox cdgd} are generated for a network with communication range $R = 0.8$, since the C-DGD algorithm does not converge for $R = 0.5$. It is remarked that as the range decreases, fewer measurements are available, making network localization more difficult.

At this point, it may be useful to also compare the complexity and run-times of the three algorithms. In general, the majorized ADMM is the most complex, since it involves solving a convex optimization problem at every iteration. On the other hand, both the C-DGD and the proximal ADMM algorithms require simple calculations with the gradients of $g_k$. For the localization problem considered here, the per-iteration majorization update at node $k$ requires solving a system of $|\mathcal{N}_k|\times |\mathcal{N}_k|)$ linear equations, incurring almost four times as much CPU-time as the proximal ADMM and the C-DGD algorithms. This extra complexity more than compensates for the reduced number of iterations afforded by the majorized ADMM. In summary, the proximal ADMM algorithm takes the least wall-time to reach a certain accuracy, say $\phi(t) = 10^{-15}$. }

\subsection{Choice of parameters}
 In order to further study the convergence rates, performance of Algorithm \ref{async_algo1} is studied for different parameter values. The convergence rate is analyzed by plotting the stopping critereon
given by $\psi(t):=\norm{\z^{t+1}-\z^t}$ against the iteration index $t$.

\begin{figure}[tbh]
			\centering
			\includegraphics[scale=.75]{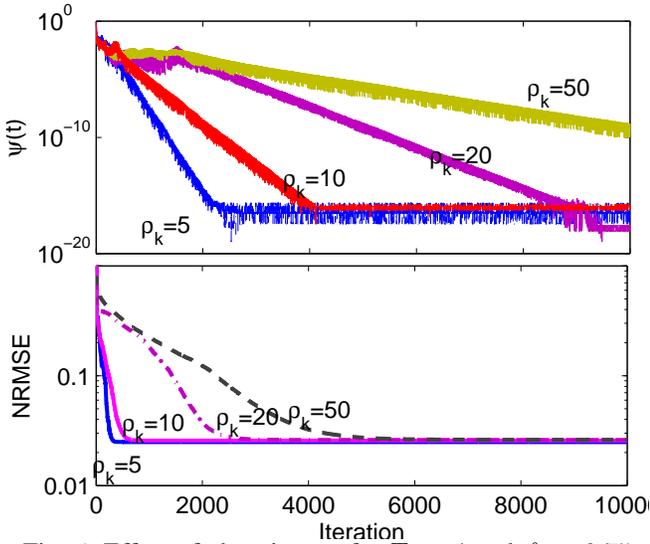}
			\caption{Effect of choosing $\rho_k$ for $T_k = 4$ and $f_k = 0.75$}
			\label{rhok_variation_combined}
\end{figure}

\subsubsection{Choice of $\rho_k$} 
Theorem \ref{thm1} guarantees the convergence of Algorithm \ref{async_algo1} whenever $\rho_k$ is chosen  in accordance with \eqref{asrho}. In practice however, it may be possible to improve the convergence rate of Algorithm \ref{async_algo1} by choosing smaller values of $\rho_k$. Fig. \ref{rhok_variation_combined} shows the evolution of $\psi(t)$ and NRMSE with iterations, for different values of $\rho_k$, while $T_k = 4$ and $f_k = 0.75$ are kept fixed. As expected, the algorithm takes longer to converge for larger values of $\rho_k$, although the asymptotic NRMSE performance for all cases remains the same. 

\begin{figure}[tbh]
	\centering
	\includegraphics[scale=.75]{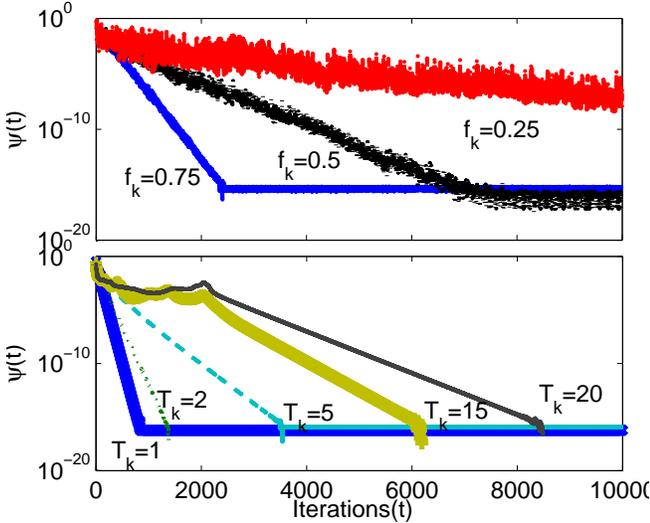}
	\caption{Effect of asynchrony for $\rho_k = 10$. The top figure uses $T_k = 4$, while the bottom figure uses $f_k = 0.75$.}
	\label{TKfkVariation}
	\end{figure}

\subsubsection{Effect of Asynchrony}
Fig. \ref{TKfkVariation} shows the effect of choosing the parameters $T_k$ and $f_k$ on the rate of convergence. To this end, we set $\rho_k = 10$ for all $k \in \cN$, and vary $T_k$ and $f_k$ separately. Both figures confirm the intuition that introduction of asynchrony results in slower convergence, even when $\rho_k$ stays the same. Note however that in order to guarantee convergence, it is necessary to choose increasingly larger values of $\rho_k$ for larger values of $T_k$ or smaller values of $f_k$. Such a choice may therefore result in even slower convergence but allow higher asynchrony.

\section{Conclusion}\label{conclusion}
This paper develops an asynchronous distributed ADMM algorithm that is applicable to a class of non-convex optimization problems. The non-convexity of the cost functions is handled either by making a first order approximation or via majorization, resulting in two variants of the proposed algorithm. Both variants converge to a stationary point of the optimization problem as long as the ADMM updates are applied ``often enough.'' The proposed algorithms find applications in distributed in-network estimation and localization. Comparisons with state-of-the-art distributed algorithms for the problem of cooperative localization in ad hoc networks demonstrates the superior performance of the proposed algorithm.

\appendices

\section{Proof of  Lemma \ref{lem1}} \label{{lem1}_proof}
This appendix provides the convergence analysis for the Lagrangian in 
Algorithm \ref{async_algo1}. Before proceeding with the proof, some 
notation is introduced. Recall the definitions of $\x_k$, $\y_k$, and 
$\z_k$, and similarly define $[\check{\x}_j]_k := x_{kj}$ and 
$[\check{\y}_j]_k := x_{kj}$ for all $k \in \N_j$. Since the Lagrangian 
is separable in both $\{\x_k\}$ and $\{z_j\}$, the following notation is 
introduced for the summands in \eqref{lagrangian} and \eqref{lang2}:
\begin{subequations}\label{lagrangian_both_fun}
	\begin{align}\label{lagrangian_fun}
	\ell_k(\x_k,\z_k,\y_k)    &:= { g_{ k }( \x_{ k }) }  +\sum _{ j 
		\in \N_k }{ \langle y_{ kj }, x_{ kj }-z_j \rangle  } \nonumber \\ 
	& \hspace{3mm}+\sum _{ j \in 
		\N_k  }{ \frac { \rho _{ k } }{ 2 } \norm{  x_{ kj }-z_j } ^{ 2 } }  \\ 
	\label{lang_kh}
	\check{\ell}_j(\check{\x}_j,z_j,\check{\y}_j)     &:= \sum _{ k 
		\in \N_j }{ \langle y_{ kj }, x_{ kj }-z_j \rangle  } \nonumber \\ 
	& \hspace{3mm}+\sum _{ k \in 
		\N_j  }{ { \frac { \rho _{ k } }{ 2 } \norm{  x_{ kj }-z_j } ^{ 2 } } } 
	+ h_j(z_j)
	\end{align}
\end{subequations}
so that $L(\{\x_k\},\z, \{\y_k\})  = \sum_{k} \ell_k(\x_k,\z_k,\y_k) + 
\sum_{k}h_k(z_k) = \sum_{j} g_j(\x_j) + \sum_{j} \check{\ell}_j 
(\check{\x}_j,z_j,\check{\y}_j)$.  Further, the approximation of 
$\ell_k(\x_k,\z_k,\y_k)$ at $\x_k = \tilde{\z}$ is denoted by [cf. 
\eqref{x_k_sync1}]
\begin{align}\label{u_k_def1}
u_k(\x_k,\tilde{\z},&\z_k,\y_k):= { g_{ k }(\tilde{\z}) }+ \langle 
\nabla g_{ k }(\tilde{\z}), \x_{ k }-\tilde{\z} \rangle\nonumber \\ & \hspace{-6mm}+ \sum _{ j \in 
	\N_k }{ \langle y_{ kj }, x_{ kj }-z_j \rangle  }+\sum _{ j\in \N_k  }{ 
	{ \frac { \rho _{ k } }{ 2 } \norm{x_{ kj }-z_j} ^{ 2 } } }.
\end{align}
With this definition, observe that the update for $\x_k^{t+1}$ in \eqref{x_k_sync1}
is given by
\begin{align} \label{xkuk}
\x_k^{t+1} &= \arg \min_{ \x_k} u_k(\x_k,\z_k^{[t+1]},\z_k^{t+1},\y_k^t) \\
\Rightarrow &\nabla_{\x_k} 
u_k(\x_k^{t+1},\z_k^{[t+1]},\z_k^{t+1},\y_k^t) = 0 \label{gradu}
\end{align}
where the gradient with respect to $\x_k$ is defined similar to that in 
\eqref{grad_syn}. In order to show that the Lagrangian decreases over 
several iterations, express the difference between consecutive 
Lagrangian values as
\begin{align}\label{Lt1}
 & L( \{\x_{ k }^{ t+1 }\},\z^{t+1},\{\y_k^{ t+1 }\})- L( \{\x_{ k }^{ 
	t }\}, \z^{t}, \{\y_k^{ t }\}) \nonumber\\
& = L( \{\x_{ k }^{ t+1 }\},\z^{t+1},\{\y_k^{ t+1 }\})-L( \{\x_{ k 
}^{t+1 }\},\z^{t+1},\{\y_k^{ t }\}) \nonumber\\
& \quad + L( \{\x_{ k }^{t+1 }\},\z^{t+1},\{\y_k^{ t }\})-L( \{\x_{ k 
}^{ t }\},\z^{t+1}, \{\y_k^{ t }\}) \nonumber\\
& \quad \quad + L( \{\x_{ k }^{t }\},\z^{t+1},\{\y_k^{ t }\})-L( \{\x_{ 
	k }^{ t }\},\z^{t}, \{\y_k^{ t }\}).
\end{align}
The subsequent lemma establishes bounds on the different terms in 
\eqref{Lt1}, and will be utilized to prove Lemma \ref{lem1}(a).
\begin{lem} \label{lemint}
	Define $\S^{t}$ as the set of nodes for which the update in 
	\eqref{zasynch} is carried out at time $t$. Then it holds that
	\begin{align}
	& L( \{\x_{ k }^{ t+1 }\},\z^{t+1},\{\y_k^{ t+1 }\})-L( \{\x_{ k }^{t+1 
	}\},\z^{t+1},\{\y_k^{ t }\})   \nonumber \\
	 &\leq \sum_{k=1}^{K}\frac { L_k^2(T_k+1) }{ 
	\rho _{ k } } \sum _{ j\in \N_k }  \sum_{d=0}^{T_k} \norm{ 
	z_{j}^{t+1-d}-z_{j}^{t-d}}^2 \label{lemy}\\
&L( \{\x_{ k }^{t+1 }\},\z^{t+1},\{\y_k^{ t }\})-L( \{\x_{ k }^{ t 
}\},\z^{t+1}, \{\y_k^{ t }\}) \nonumber\\
&\quad\leq -\sum_{k=1}^K{ \frac { \rho _{ k }-7L_k }{ 2 } \sum_{j\in 
		\N_k}\|x_{kj}^{t+1}-x_{kj}^{t}\|^2 }\nonumber \\
&\hspace{-3mm} + \sum_{k=1}^{K}\left(\frac { 7 
	L^3_k(T_k+1)}{ 2 \rho _{ k }^2 
}+\frac{L_kT_k}{2}\right)\sum_{d=0}^{T_k-1}\sum_{j\in 
\N_k}\|z_{j}^{t+1-d}-z_{j}^{t-d}\|^2 \label{lemx}\\
& L( \{\x_{ k }^{t }\},\z^{t+1},\{\y_k^{ t }\})-L( \{\x_{ k }^{ t 
}\},\z^{t}, \{\y_k^{ t }\})  \nonumber \\
&\leq -\sum_{k \in \S^{t}} \frac{\rho_k}{2} 
\norm{z_k^{t+1}-z_k^t}^2 \label{lemz}
\end{align}
\end{lem}

\begin{IEEEproof} The proof begins by establishing a bound on difference 
	between successive dual values. Rearranging \eqref{xasynch},
	\begin{align}
	[\nabla g_k(\z_k^{[t+1]})]_j + y_{kj}^t + \rho_{k}( 
	x_{kj}^{t+1}-z_j^{t+1})&= 0, & j &\in \N_k \nonumber
	\end{align}
	which, together with the update for $y_{kj}^{t+1}$ [cf. \eqref{eq_y_k1}
	] yields
	\begin{align}\label{lemma3.1a}
	[\nabla g_k({\z}_k^{[t+1]})]_j&=-y_{kj}^{t+1}, & j &\in \N_k.
	\end{align}
	Similarly, it holds that $y_{kj}^t = -[\nabla 
	g_k({\z}_k^{[t]})]_j$ for all $j \in \N_k$. Therefore, for each $k 
	\in \cK$, the following bound applies:
	\begin{subequations}\label{ybnd}
		\begin{align}
		\sum_{j \in \N_k} & \norm{y_{kj}^{t+1}-  y_{kj}^{t}}^2 =    \norm{ \nabla 
			g_k(\z_k^{[t+1]}) -\nabla g_k(\z_k^{[t]})}^2 \nonumber\\
		& \leq L_k^2 \norm{\z_k^{[t+1]}- \z_k^{[t]}}^2 \label{lip1}\\
		& \leq L_k^2 \left(\sum_{d=0}^{T_k} \norm{\z_k^{t+1-d}- \z_k^{t-d}} 
		\right)^2 \label{triangle1}\\
		& \leq L_k^2(T_k+1) \sum_{d=0}^{T_k} \sum_{j\in 
			\N_k}\norm{z_{j}^{t+1-d}-z_{j}^{t-d}}^2 \label{triangle2}
		\end{align}
	\end{subequations}
	where \eqref{lip1} follows from Assumption (A1) while \eqref{triangle1} 
	and \eqref{triangle2} follow from the use of triangle inequality. Next, 
	the different bounds in \eqref{lemy}, \eqref{lemx}, and \eqref{lemz} are 
	proved.
	
	\indent (a) The bound in \eqref{lemy} follows from the following 
	equalities utilizing \eqref{eq_y_k1}.
	\begin{align}
&	L( \{\x_{ k }^{ t+1 },\z^{t+1},\{\y_k^{ t+1 }\})  - L( \{\x_{ k }^{t+1 
	}\},\z^{t+1},\{\y_k^{t}\}) \nonumber\\
	& =\sum _{ k=1 }^{ K }\sum _{ j\in\N_k }  \langle y_{ kj }^{ t+1 }, x_{ 
		kj }^{ t+1 }-z_j^{t+1} \rangle  - \langle y_{ kj }^{ t }, x_{ kj }^{ t+1 
	}-z_j^{t+1} \rangle
	\nonumber\\
	& =\sum _{ k=1 }^{ K }\sum _{ j\in\N_k }{ \langle y_{ kj }^{ t+1 }-y_{ 
			kj }^{ t },\frac{1}{\rho_{k}}\{y_{ kj }^{ t+1 }-y_{ kj }^{ t }\} \rangle  }
	\nonumber\\
	& =\sum_{k=1}^{K}  { \frac { 1 }{ \rho _{ k } }  \sum _{ j\in\N_k 
		}\left\| y_{ kj }^{t+1}-y^{ t }_{kj} \right\| ^{ 2 } } \label{lemyproof}
	\end{align}
	Finally, the bound in \eqref{lemy} follows by substituting 
	\eqref{triangle2} into the right-hand side of \eqref{lemyproof}.
	
	\indent (b) First some properties of the individual summands 
	$\ell_k(\x_k,\z_k^{t+1},\y_k^t)$ and their approximated versions 
	$u_k(\x_k, \tilde{\z}, \z_k^{t+1}, \y_k^t)$ are established. Since 
	$\nabla g_k(\x_k)$ is Lipschitz continuous in $\x_k$, the following 
	quadratic upper bound holds for all $\x_k$
	\begin{align}
	g_k( \x_k) \leq g_k(\z_k^{t+1})&+ \langle \nabla g_k(\z_k^{t+1}) , \x_{ k 
	}-\z_k^{t+1} \rangle \nonumber \\ &+\frac{L_k}{2}\sum_{j \in 
	\N_k}\norm{x_{kj}-z_j^{t+1}}^{ 2 } \label{lipcong}
\end{align}
which implies that
\begin{align}\label{lemma2_1}
\ell_k(\x_k,{\z}_k^{t+1},\y_k^t) \le 
u_k(\x_k,&\z_k^{t+1},\z_k^{t+1},\y_k^t)\nonumber \\ & + \frac{L_k}{2}\sum_{j \in 
	\N_k}\norm{x_{kj}-z_j^{t+1}}^{ 2 }.
\end{align}

Observing that $u_k(\x_k,\tilde{\z},\z_k^{t+1},\y_k^t)$ is 
differentiable and strongly convex with respect to $\x_k$, it holds that
\begin{align}
& u_k(\x_k^{t+1},\tilde{\z},\z_k^{t+1},\y_k^{t}) - 
u_k(\x_k^{t},\tilde{\z},\z_k^{t+1},\y_k^{t}) \nonumber\\
& \leq \langle \nabla_{\x_k} 
u_k(\x_k^{t+1},\tilde{\z},\z_k^{t+1},\y_k^{t}),x_{kj}^{t+1}-x_{kj}^t 
\rangle \nonumber \\ & -\frac{\rho_k    }{2}\norm{\x_k^{t+1}-\x_k^{t}}^2 \label{scu}
\end{align}
for all $\tilde{\z}$. In particular for $\tilde{\z} = \z_k^{[t+1]}$, the 
gradient on the right is zero from \eqref{gradu}, so that
$ u_k(\x_k^{t+1},\z^{[t+1]},\z_k^{t+1},\y_k^{t}) - u_k(\x_k^{t},\z^{[t+1]},\z_k^{t+1},\y_k^{t})
\leq  -\frac{\rho_k}{2}\sum_{j \in \N_k} \norm{x_{kj}^{t+1}-x_{kj}^t}. $
Note also that
$ \nabla_{\x_k} 
u_k(\x_k^{t+1},\z_k^{[t+1]},\z_k^{t+1},\y_k^{t}) -\nabla_{\x_k} 
u_k(\x_k^{t+1},\z_k^{t+1},\z_k^{t+1},\y_k^{t}) = \nabla g_{ k }(\z_{ k 
}^{[t+1]})-\nabla g_{ k }( \z_{ k }^{t+1}) $

Next, specializing \eqref{scu} for $\tilde{\z} = \z_k^{t+1}$, the 
following series of inequalities are obtained
\begin{subequations}\label{lemma3_main}
	\begin{align}
	u_{ k }&(\x_k^{t+1},\z_k^{t+1},\z_k^{t+1},\y_k^t) - u_{ k 
	}(\x_k^{t},\z_k^{t+1},\z_k^{t+1},\y_k^t)  \nonumber \\
	& \leq  \langle \nabla_{\x_k} u_{ k 
	}(\x_k^{t+1},\z_k^{t+1},\z_k^{t+1},\y_k^t) ,  \x^{ t+1 }_k-\x^{t}_{ k } 
	\rangle \nonumber \\ &\qquad -  \frac { \rho _{ k } }{ 2 } \norm{\x_{ k }^{t+1}- \x^{ t 
		}_k}^{ 2 }  \label{lemma3_2eqn2_a} \\
	& = {\langle  \nabla_{\x_k} u_{ k 
		}(\x_k^{t+1},\z_k^{t+1},\z_k^{t+1},\y_k^t),  \x^{ t+1 }_k-\x^{t}_{ k 
} \rangle}\nonumber \\ 
&- {\langle \nabla_{\x_k} u_{ k 
}(\x_k^{[t+1]},\z_k^{t+1},\z_k^{t+1},\y_k^t) ,  \x^{ t+1 }_k-\x^{t}_{ k 
} \rangle} \nonumber \\ 
&\qquad-{ \frac { \rho _{ k } }{ 2 } \left\| \x_{ k }^{t+1}- \x^{ t 
}_k \right\| ^{ 2 } } \label{lemma3_2eqn2_b}\\
& = \langle \nabla g_{ k }( \z_{ k }^{t+1})-\nabla g_{ k }( \z_{ 
	k }^{[t+1]}),  \x_{ k }^{t+1}- \x^{ t }_k \rangle \nonumber \\ 
& -\qquad { \frac { \rho _{ k } 
	}{ 2 } \norm{ \x_{ k }^{t+1}- \x^{ t }_k }^{ 2 } } \nonumber \\
&\leq L_k\norm{ \z_{ k }^{[t+1]}-\z_{ k }^{t+1}}\norm{\x_{ k 
	}^{t+1}- \x^{ t }_k}-{ \frac { \rho _{ k } }{ 2 } \norm{\x_{ k }^{t+1}- 
	\x^{ t }_k} ^{ 2 } } \label{lemma3_2eqn2_c} \\
&\leq \frac{L_kT_k}{2}\sum_{d=0}^{T_k-1}\norm{\z_{ k }^{t+1-d}-\z_{ 
		k }^{t-d}}^2 -  \frac{ \rho _{ k } -L_k}{ 2 } \norm{\x_{ k }^{t+1}- \x^{ 
		t }_k} ^{ 2 } \nonumber \\
& = \frac{L_kT_k}{2}\sum_{d=0}^{T_k-1}\sum_{j\in 
	\N_k}\|z_{j}^{t+1-d}-z_{j}^{t-d}\|^2\nonumber \\ 
&\qquad -{ \frac { \rho _{ k } -L_k}{ 2 } 
	\sum_{j\in \N_k} \|  x_{ kj }^{t+1}- x^{ t }_{kj} \| ^{ 2 } } 
\label{lemma3_2eqn2_e}
\end{align}
\end{subequations}
where \eqref{lemma3_2eqn2_b} follows from \eqref{gradu}, 
\eqref{lemma3_2eqn2_c} follows from (A1), and the rest follow from the 
use of triangle inequality similar to that in \eqref{ybnd}.

Similarly, observe that
\begin{align}\label{lemma3_subseq_1}
  & u_{ k }(\x_k^{t+1},\z_k^{t+1},\z_k^{t+1},\y_k^t)-\ell_{ k 
}(\x_k^{t+1},\z_k^{t+1},\y_k^t)\nonumber \\
 &=g_k(\z_k^{t+1})-g_k( 
\x_k^{t})+{ \langle \nabla g_k(\z_k^{t+1}) , \x_k^{t}-\z_k^{t+1} \rangle}.
\end{align}
Similar to \eqref{lipcong}, the following quadratic upper bound is also 
implied by (A1),
\begin{align}
 g_k(\z_k^{t+1}) \leq g_k(\x_k^t) + & \langle \nabla g_k(\x_k^t) , 
\z_k^{t+1} - \x^t_{ k }\rangle \nonumber \\ 
&+\frac{L_k}{2}\sum_{j \in 
	\N_k}\norm{x_{kj}^t-z_j^{t+1}}^{ 2 }. \label{lipcong2}
\end{align}
Substituting \eqref{lipcong2} into \eqref{lemma3_subseq_1}, we obtain
\begin{subequations}\label{lem3.21b}
	\begin{align}
	& u_{ k }(\x_k^{t+1},\z_k^{t+1},\z_k^{t+1},\y_k^t)-\ell_{ k 
	}(\x_k^{t+1},\z_k^{t+1},\y_k^t) \nonumber \\ & \leq  \langle \nabla 
	g_k({\z}_k^{t+1})-\nabla g_k( \x_k^{t}) , \x_k^{t}-{\z}_k^{t+1} \rangle 
 \nonumber \\ &	+\frac{L_k}{2}\sum_{j \in \N_k}\norm{x_{kj}^t-z_j^{t+1}}^{ 2 }\nonumber\\
	&  \le \frac{3L_k}{2}\sum_{j \in \N_k}\norm{x_{kj}^t-z_j^{t+1}}^{ 2 } 
	\label{lem3.21b1}\\
	& \le 3 L_k \sum_{j \in \N_k}\left(\norm{x_{kj}^{t}- x_{ kj }^{ t+1 }}^{ 
		2 }+\norm{x_{kj}^{t+1}-z_j^{t+1}}^{ 2 }\right)\label{lem3.21b2}
	\end{align}
\end{subequations}
where \eqref{lem3.21b1} follows from (A1) and \eqref{lem3.21b2} from the 
use of triangle inequality.

Having derived inequalities \eqref{lemma2_1}-\eqref{lem3.21b}, the 
difference between the summands of \eqref{lemx} can finally be bounded. 
Towards this end, it holds that
\begin{align}
&\ell_{ k }( \x_{ k }^{ t+1 },\z_k^{t+1},\y_k^{ t })-\ell_{ k 
}(\x_k^{t},\z_k^{t+1},\y_k^{ t })\nonumber\\ & \le u_k( 
\x_k^{t+1},\z_k^{t+1},\z_k^{t+1},\y_k^t) - u_k( 
\x_k^{t},\z_k^{t+1},\z_k^{t+1},\y_k^t) \nonumber \\ & + \frac{L_k}{2}\sum_{j\in \N_k} 
\norm{z_j^{t+1}- x^{t+1}_{kj}}^{ 2 } \nonumber\\
& \quad + u_k( \x_k^{t},\z_k^{t+1},\z_k^{t+1},\y^t)-\ell_k( 
\x_k^{t},\z_k^{t+1},\y^t)\nonumber\\
&\leq  \frac{L_kT_k}{2}\sum_{d=0}^{T_k-1}\sum_{j\in 
	\N_k}\|z_{j}^{t+1-d}-z_{j}^{t-d}\|^2 \nonumber \\ 
&-{ \frac { \rho _{ k } -L_k}{ 2 } 
	\sum_{j\in \N_k} \|  x_{ kj }^{t+1}- x^{ t }_{kj} \| ^{ 2 } }+ 
\frac{L_k}{2}\sum_{j\in \N_k} \norm{z_j^{t+1}- x^{t+1}_{kj}}^{ 2 }  \nonumber\\
&  + 3 L_k \sum_{j \in \N_k}\left(\norm{x_{kj}^{t}- x_{ kj }^{ t+1 
	}}^{ 2 }+\norm{x_{kj}^{t+1}-z_j^{t+1}}^{ 2 }\right) \\
	&= -{ \frac { \rho _{ k }-7L_k }{ 2 } \sum_{j\in 
			\N_k}\|x_{kj}^{t+1}-x_{kj}^{t}\|^2 } +\frac { 7 L_k }{ 2 \rho _{ k }^2 } 
	\sum_{j\in 
		\N_k}\|y_{kj}^{t+1}-y_{kj}^{t}\|^2
\nonumber \\ 
&\qquad+\frac{L_kT_k}{2}\sum_{d=0}^{T_k-1}\sum_{j\in 
		\N_k}\|z_{j}^{t+1-d}-z_{j}^{t-d}\|^2 \label{ysub}\\
	&= -{ \frac { \rho _{ k }-7L_k }{ 2 } \sum_{j\in 
			\N_k}\|x_{kj}^{t+1}-x_{kj}^{t}\|^2 } \nonumber \\ 
	&\qquad +\left(\frac { 7 L^3_k(T_k+1)}{ 2 
		\rho _{ k }^2 }+\frac{L_kT_k}{2}\right)\sum_{d=0}^{T_k-1}\sum_{j\in 
		\N_k}\|z_{j}^{t+1-d}-z_{j}^{t-d}\|^2 \label{ysub2}
	\end{align}
	where \eqref{ysub} follows from using the update for $y_{kj}^{t+1}$ and 
	\eqref{ysub2} from \eqref{ybnd}. Finally, summing \eqref{ysub2} over $k 
	= 1, 2, \ldots, K$, the result in \eqref{lemx} follows.
	
	\indent (c) Define the indicator function $\ind_{\cX}(x) = 1$ if $x \in 
	\cX$ and zero otherwise, the observe that update for $z_j$ can be 
	written as [cf. \eqref{z_k_update}],
	\begin{align}
	z_j^{t+1} &= \arg\min_{z_j \in \cX_j} \check{\ell}(\x_j^{t},z_j,\y_j^t) \label{zind0}\\
	&= \arg\min_{z_j} \check{\ell}(\x_j^{t},z_j,\y_j^t) + \ind_{\cX_j}(z_j)  
	\label{zind}
	\end{align}
	for all $j \in \S^t$. The first order optimality condition for 
	\eqref{zind} is that
	\begin{align}
	0 \in \partial \left(\check{\ell}(\x_j^{t},z^{t+1}_j,\y_j^t) + 
	\ind_{\cX_j}(z^{t+1}_j)\right)\label{fooc}
	\end{align}
	where $\partial(\cdot)$ denotes the subgradient operator. Observe 
	further that since the function $\check{\ell}(\x_j^{t},z_j,\y_j^t) + 
	\ind_{\cX_j}(z_j)$ is strongly convex in $z_j$, it holds that
	\begin{align}
	& \check{\ell}(\x_j^{t},z^{t+1}_j,\y_j^t) + \ind_{\cX_j}(z^{t+1}_j) 
	\nonumber\\
	& \qquad \leq \check{\ell}(\x_j^{t},z^{t}_j,\y_j^t) + 
	\ind_{\cX_j}(z^{t}_j) - 
	\frac{\rho_j}{2}\norm{z_{j}^{t+1}-z_j^t}^2\nonumber \\
	& - \langle \partial 
	\left(\check{\ell}(\x_j^{t},z^{t+1}_j,\y_j^t) + 
	\ind_{\cX_j}(z^{t+1}_j)\right), z_j^t-z_j^{t+1}\rangle. \label{strconv1}
	\end{align}
	Since $\ind_{\cX_j}(z_j^{t+1}) = \ind_{\cX_j}(z_j^t) = 0$ for $ j \in \S^t$, \newline  from 
	\eqref{fooc} and \eqref{strconv1} it holds that
	\begin{align}
	\check{\ell}_j(\check{\x}_j,z_j^{t+1},\check{\y}_j) - 
	\check{\ell}_j(\check{\x}_j,z^{t}_j,\check{\y}_j)  \leq 
	-\frac{\rho_j}{2}\norm{z_{j}^{t+1}-z_j^t}^2  \label{lcj}
	\end{align}
	For all $j \notin \S^t$, since $z_j^{t+1}=z_j^t$, the right-hand side is 
	clearly zero. Finally, the result in \eqref{lemz} is obtained by summing 
	both sides in \eqref{lcj} over $j = 1, 2, \ldots, K$.
	\end{IEEEproof}
	
	\begin{IEEEproof}[Proof of Lemma \ref{lem1}(a)]
		From Lemma \ref{lemint}, the decrease in the Lagrangian over 
		consecutive time slots is given by
		\begin{align} \label{consecdiff}
		&L( \{\x_{ k }^{ t+1 }\},{\z}^{t+1},\{\y_k^{ t+1 }\})-L( \{\x_{ k }^{ t 
		}\},\z^{t}, \{\y_k^{ t }\}) \nonumber\\
		& \leq -\sum_{k=1}^{K}{ \frac { \rho _{ k }-7L_k }{ 2 } \sum_{j \in 
				\N_k} \norm{x_{kj}^{t+1}- x_{kj}^t}^2 } - \sum_{k\in 
			\set^{t}}\frac{\rho_k}{2} \norm{z_{k}^{t+1}-z_{k}^{t}}^2 \nonumber\\
		& + \sum_{k=1}^{K} \left(\frac{1    }{\rho_k}+\frac { 7L_k }{ 2\rho_k^2 
		}\right)L_k^2(T_k+1) \sum_{d=0}^{T_k} \sum_{j\in 
		\N_k}\norm{z_{j}^{t+1-d}-z_{j}^{t-d}}^2 \nonumber \\ &\qquad + \sum_{k=1}^K 
	\frac{L_kT_k}{2}   \sum_{d=0}^{T_k-1} \sum_{j\in 
		\N_k}\norm{z_{j}^{t+1-d}-z_{j}^{t-d}}^2
	\end{align}
	Given $t_0 \geq 1$ and $T \geq 1$, define the inverse mapping $\S^{-1}_T(k) := \{ t_0 \leq t\leq t_0 + T \mid k \in \S^t \}$ as the set of iterations for which update \eqref{z_k_update} is 
	applied at node $k$ and let $f_k \geq \abs{\S^{-1}_T(k)}/T$. Then, summing both sides of \eqref{consecdiff} over $t = t_0$, $t_0 + 1$, $\ldots$, $T+t_0$,
	\begin{align} \label{convfint}
&	L( \{\x_{ k }^{ T+ t_0 }\},{\z}^{T+t_0},\{\y_k^{T+t_0}\})  -L( \{\x_{ k 
	}^{ t_0 }\},{\z}^{t_0}, \{\y_k^{t_0}\}) \leq\nonumber \\
 -	& \hspace{-2mm}\sum_{t=t_0}^{T+t_0}\sum_{k=1}^K { \frac{\beta_k	}{2} 
		\sum_{j\in \N_k}\norm{x_{kj}^{i+1}-x_{kj}^{i}}^2 } + 
	\alpha_k\norm{z_{k}^{i+1}-z_{k}^{i}}^2 
	\end{align}
	where $\alpha_k $ and $\beta_k$ are defined as in \eqref{rholk}, yielding the desired result. Note that from Assumption (A4), the term on 
	the right-hand side of \eqref{convfint} is negative.

\end{IEEEproof}

\begin{IEEEproof}[Proof of Lemma 1(b)]: Using the Lipschitz continuity 
	of $\nabla g_k(\cdot)$ and applying triangle inequality, it follows that
	\begin{align}\label{gkbnd}
&	g_k(\z_k^{t+1}) \leq  g_{ k }( \x_{ k }^{t+1}) + \langle\nabla g_{ k }( 
	\x_k^{t+1}),\z_k^{t+1}- \x_{ k }^{t+1} \rangle \nonumber \\ & \qquad + \frac { L_{ k } }{ 2 } 
	\norm{\x_{ k }^{t+1}-{\z}_k^{t+1}}^{ 2 }\\
	&= g_{ k }( \x_{ k }^{t+1})+\langle\nabla g_{ k }( \x_k^{t+1})-\nabla 
	g_{ k }(\z_k^{t+1}),\z_k^{t+1}- \x_{ k }^{t+1} \rangle \nonumber\\
	& + \langle\nabla g_{ k }(\z_k^{t+1}),\z_k^{t+1}- \x_{ k }^{t+1} 
	\rangle  + \frac { L_{ k } }{ 2 } \norm{\x_{ k }^{t+1}-\z_k^{t+1}}^{ 2 } \\
	& \leq  g_{ k }( \x_{ k }^{t+1})+\langle\nabla g_{ k 
		}({\z}_k^{t+1}),{\z}_k^{t+1}- \x_{ k }^{t+1} \rangle \nonumber \\ & \qquad+ \frac { 3L_{ k } 
	}{ 2 } \norm{\x_{ k }^{t+1}-{\z}_k^{t+1}} ^{ 2 }
\end{align}

Next, using the relationship from \eqref{lemma3.1a}, observe that the 
\begin{align}\label{l34a}
& L(\{ \x_{ k }^{t+1}\}, {\z}^{t+1};\{\y_k^{ t+1 }\})\nonumber\\
& =\sum_{ k=1 }^{ K }\Bigg({ g_{ k }( \x_{ k }^{t+1}) }  +\sum _{ j \in \N_k }{ 
	\langle y_{ kj }^{t+1}, x_{ kj }^{t+1}-z_j^{t+1} \rangle  }\nonumber \\ &\qquad  + \frac { 
	\rho _{ k } }{ 2 } \sum _{ j \in \N_k  }  \norm{ x_{ kj 
	}^{t+1}-z_j^{t+1} } ^{ 2 } +h_k(z_k^{t+1})\Bigg) \\
&   \ge\sum _{ k=1 }^{ K } \Bigg( g_{ k 
	}(\z_k^{t+1})+h_k(z_k^{t+1})\nonumber \\ &\qquad+\frac { \rho _{ k }-3L_k }{ 2 }\sum _{ j 
	\in \N_k  }  \norm{ x_{ kj }^{t+1}-z_j^{t+1} } ^{ 2 }\nonumber \\ &\qquad  + \langle\nabla 
g_{ k }({\z}_k^{t+1}) - \nabla g_{ k }({\z}_k^{[t+1]}),{\z}_k^{t+1}- 
\x_{ k }^{t+1} \rangle \Bigg)  
\end{align}
Next, using Cauchy-Schwarz inequality on the last term, it follows that
\begin{align}
& L(\{ \x_{ k }^{t+1}\}, {\z}^{t+1};\{\y_k^{ t+1 }\})\nonumber \\ & \geq \mathsf{P} + 
\frac { \rho _{ k }-3L_k }{ 2 }\sum _{ j \in \N_k  } \norm{ x_{ kj 
	}^{t+1}-z_j^{t+1} } ^{ 2 } \nonumber \\
&\qquad - \norm{\nabla g_{ k }({\z}_k^{t+1}) - \nabla 
	g_{ k }({\z}_k^{[t+1]})}\norm{{\z}_k^{t+1}- \x_{ k }^{t+1}} \\
&\geq \mathsf{P} + \frac { \rho _{ k }-4L_k }{ 2 }\sum _{ j \in \N_k  }  
\norm{ x_{ kj }^{t+1}-z_j^{t+1} } ^{ 2 } \nonumber \\
&\qquad- 
\frac{L_k}{2}\sum_{j\in\N_k}\norm{z_j^{t+1} - z_j^{[t+1]}}^2 
\label{lbound2}\\
& \geq \mathsf{P} -\frac{L_k}{2}\sum_{j\in \N_k}\text{diam}^2(\cX_j) > -\infty
\end{align}
where the last inequality follows from the fact that $\rho_k \geq 7L_k$ 
[cf. (A4)], and that $\cX_j$ is compact [cf. (A2)].
\end{IEEEproof}

\section{Proof of  Theorem \ref{thm1}} \label{{theorem}_proof}
\begin{IEEEproof}[Proof of Theorem \ref{thm1}(a)]
From Lemma \ref{lem1}, it follows that $L(\{ \x_{ k }^{t}\}, {\z}^{t};\{\y_k^{ t}\})$ converges as $t \rightarrow \infty$. Therefore, it holds from \eqref{lagdiff} that,
\begin{subequations} \label{thmproof1}
\begin{align}
\lim_{ t\rightarrow \infty  }{\left\| x_{kj}^{t+1} -  x_{kj}^{t} \right\|  } &\rightarrow 0 & j\in \N_k, & \quad k = 1,\ldots, K \label{thmproof1a}\\ 
\lim_{ t\rightarrow \infty  }{\left\| z_k^{t+1} -  z_k^{t} \right\|  } & \rightarrow 0 && k = 1, \ldots, K \label{thmproof1b} 
\end{align}
Next, using \eqref{thmproof1b} into \eqref{ybnd} and the form of the update in \eqref{eq_y_k1}, we obtain
\begin{align}
\lim_{ t\rightarrow \infty  }{\left\| y_{kj}^{t+1} -  y_{kj}^{t} \right\|  }\rightarrow & 0 &j\in \N_k, \quad k = 1, \ldots, K \label{thmproof1c}\\ 
\lim_{ t\rightarrow \infty  }{\left\| x_{kj}^{t} -  z_{j}^{t} \right\|  }\rightarrow & 0 & j\in \N_k, \quad k =  1, \ldots, K. \label{thmproof1d}
\end{align} 
\end{subequations}
\end{IEEEproof}

\begin{IEEEproof}[Proof of Theorem \ref{thm1}(b)] From \eqref{thmproof1}, it follows that the limit points $\{\{x_{kj}^\star\}, \{z_j^\star\}, \{y_{kj}^\star\}\}$ exist, and satisfy
\begin{align}
x_{kj}^\star &= z_j^\star & j \in \N_k, &~ k = 1, \ldots, K \label{primfeas}
\end{align}
which is the primal feasibility condition for \eqref{prob_formulation}. Since $z_j^{t+1} \in \cX_j$ for all $t$, it should hold that $z_j^\star \in \cX_j$ and $x_{kj}^\star=x_{lj}^\star$ for all $j,l,k\in \N_k$. The first order optimality condition for $\x_k$ can be obtained from \eqref{lemma3.1a}, which implies \eqref{th2_a}.
Finally, the first order optimality condition for $z_j$ can be obtained from \eqref{fooc}, which implies that
\begin{align}
0 \in \partial \left(h_j(z) + \ind_{z\in \cX_j}\right)\mid_{z = z_j^\star} - \sum_{k \in \N_j} y_{kj}^\star - \sum_{k\in \N_j}\rho_k\left(x_{kj}^\star - z_j^\star\right). \nonumber
\end{align} 
The desired result follows by using \eqref{primfeas} and noting that $0 \in \partial  \ind_{z \in \cX_j}\mid_{z = z_j^\star}$. 
\end{IEEEproof}

\section{Convergence of the Majorized ADMM} \label{madmm}
The proof of convergence of Theorem \ref{thm2} follows the same general outline as that of Theorem \ref{thm1}. The definitions of  $\ell_k$ and $\check{\ell}_j$ are the same as in \eqref{lagrangian_both_fun}, while the approximation of $\ell_k(\x_k,\y_k,\z_k)$ at $\x_k = \tilde{\z}$ is now defined as
\begin{align}\label{u_k_def2}
& u_k(\x_k,\tilde{\z},\z_k,\y_k):=  f_k(\x_k,\tilde{\z}) + \sum _{ j \in 
	\N_k }{ \langle y_{ kj }, x_{ kj }-z_j \rangle  }\nonumber \\ & \qquad \qquad+\sum _{ j\in \N_k  }{ 
	{ \frac { \rho _{ k } }{ 2 } \norm{x_{ kj }-z_j} ^{ 2 } } }.
\end{align}
From the update in \eqref{mxasynch}, it holds that  $\nabla_{\x_k} 
u_k(\x_k^{t+1},\z_k^{[t+1]},\z_k^{t+1},\y_k^t) = 0$, and  
\begin{align}\label{grad2}
[\nabla f_k(\x_k^{t+1},\z_k^{[t+1]})]_j + y^{t}_{kj} + \rho_k (x_kj - z_j^{t+1}) = 0
\end{align}
Further, from the update of $y_{kj}^{t+1}$ and from \eqref{grad2}, 
	\begin{align}\label{ytplus1}
	y_{kj}^{t+1} &= -[\nabla_{\x_k} f_k(\x_k^{t+1},\z_k^{[t+1]})]_j, & j &\in \N_k.
	\end{align}
	Similarly, it holds that $y_{kj}^t = -[\nabla_{\x_k} f_k(\x_k^{t},\z_k^{[t]})]_j$ for all $j \in \N_k$. Therefore, for each $k \in \cK$, the following bound applies:
	\begin{subequations}\label{ybnd_maj}
		\begin{align}
	&	\sum_{j \in \N_k}\norm{y_{kj}^{t+1}-y_{kj}^{t}}^2 =    \norm{\nabla f_k(\x_k^{t+1},\z_k^{[t+1]}) -\nabla f_k(\x_k^{t},\z_k^{[t]})}^2 \nonumber\\
		& \leq 2L_{k}^2\norm{\x_k^{t+1}-\x_k^t}^2+ 2L_k^2 \norm{\z_k^{[t+1]}- \z_k^{[t]}}^2 \label{lip21}\\
		& \leq 2L_k^2\norm{\x_k^{t+1}-\x_k^t}^2+2L_k^2(T_k+1) \sum_{d=0}^{T_k} \sum_{j\in 
			\N_k}\norm{z_{j}^{t+1-d}-z_{j}^{t-d}}^2 \label{triangle3}
		\end{align}
	\end{subequations}
	where \eqref{lip21} follows from Assumption, and \eqref{triangle3} follows similarly as in \eqref{triangle2}.

As in Appendix \ref{{lem1}_proof}, the Lagriangian is split into three summands [cf. Lemma \ref{lemint}], each of which must be separately bounded. The bound on the first summand follows directly from \eqref{ybnd_maj}. 
	\begin{align}
 &	L( \{\x_{ k }^{ t+1 },\z^{t+1},\{\y_k^{ t+1 }\})  - L( \{\x_{ k }^{t+1 
	}\},\z^{t+1},\{\y_k^{t}\}) \nonumber \\
	& =\sum_{k=1}^{K}  { \frac { 1 }{ \rho _{ k } }  \sum _{ j\in\N_k 
		}\left\| y_{ kj }^{t+1}-y^{ t }_{kj} \right\| ^{ 2 } } \\
		& \leq\sum_{k=1}^{K} \frac{2L_k^2}{\rho_k}\sum_{j\in \N_k}\norm{x_{kj}^{t+1}-x_{kj}^t}^2\nonumber \\ 
		&\hspace{2mm} +\sum_{k=1}^K\frac{2L_k^2(T_k+1)}{\rho_k} \sum_{d=0}^{T_k} \sum_{j\in 
			\N_k}\norm{z_{j}^{t+1-d}-z_{j}^{t-d}}^2.
	\end{align}
	
For the second summand, some relevant inequalities are first stated. From the definition of the majorizing function, it holds that for all $\x_k$, $\ell_k(\x_k,\z_k^{t+1},\y_k^t) \le u_k(\x_k,\z_k^{t+1},\z_k^{t+1},\y_k^t)$.
Since $f_k(\x_k,\z_k^{t+1})$ is differentiable and convex with respect to $\x_k$, it follows that $u_k(\x_k,\tilde{\z},\z_k^{t+1},\y_k^t)$ is differentiable and strongly convex, implying that

\begin{align}
 & u_k(\x_k^{t+1},\tilde{\z},\z_k^{t+1},\y_k^{t})- 
u_k(\x_k^{t},\tilde{\z},\z_k^{t+1},\y_k^{t}) \nonumber\\
& \leq \langle \nabla_{\x_k} 
u_k(\x_k^{t+1},\tilde{\z},\z_k^{t+1},\y_k^{t}),x_{kj}^{t+1}-x_{kj}^t 
\rangle\nonumber\\
&\qquad  -\frac{\rho_k}{2}\norm{\x_k^{t+1}-\x_k^{t}}^2 \label{scu2}
\end{align}
for all $\tilde{\z}$. In particular for $\tilde{\z} = \z_k^{[t+1]}$, the 
gradient on the right is zero from \eqref{grad2}, so that
$ u_k(\x_k^{t+1},\z^{[t+1]},\z_k^{t+1},\y_k^{t})  - 
u_k(\x_k^{t},\z^{[t+1]},\z_k^{t+1},\y_k^{t}) \leq 
-\frac{\rho_k}{2}\sum_{j \in \N_k} \norm{x_{kj}^{t+1}-x_{kj}^t}.
$

Note also that $  u_k(\x_k^{t+1},\z_k^{[t+1]},\z_k^{t+1},\y_k^{t})-\nabla_{\x_k} 
u_k(\x_k^{t+1},\z_k^{t+1},\z_k^{t+1},\y_k^{t})
 = \nabla_{\x_k} f_k(\x_k^{t+1},\z_k^{[t+1]}) - \nabla_{\x_k} f_k(\x_k^{t+1},\z_k^{t+1})$. Next, specializing \eqref{scu2} for $\tilde{\z} = \z_k^{t+1}$, the 
following series of inequalities are obtained
\begin{subequations}\label{list}
	\begin{align}
	u_{ k }&(\x_k^{t+1},\z_k^{t+1},\z_k^{t+1},\y_k^t) - u_{ k 
	}(\x_k^{t},\z_k^{t+1},\z_k^{t+1},\y_k^t)  \nonumber \\
	& \leq  \langle \nabla_{\x_k} u_{ k 
	}(\x_k^{t+1},\z_k^{t+1},\z_k^{t+1},\y_k^t) ,  \x^{ t+1 }_k-\x^{t}_{ k } 
	\rangle   \nonumber \\ & \qquad -\frac { \rho _{ k } }{ 2 } \norm{\x_{ k }^{t+1}- \x^{ t 
		}_k}^{ 2 }  \label{list1} \\
& \leq \frac{L_kT_k}{2}\sum_{d=0}^{T_k-1}\sum_{j\in 
	\N_k}\|z_{j}^{t+1-d}-z_{j}^{t-d}\|^2\nonumber \\ &\qquad-{ \frac { \rho _{ k } -L_k}{ 2 } 
	\sum_{j\in \N_k} \|  x_{ kj }^{t+1}- x^{ t }_{kj} \| ^{ 2 } } 
\label{list5}
\end{align}
\end{subequations}
 Utilizing \eqref{grad2} and following steps \eqref{lemma3_2eqn2_a}-\eqref{lemma3_2eqn2_e},  \eqref{list5} is obtained, which follow similarly as in \eqref{lemma3_main}. 

The Lipschitz continuity of $\nabla g_k(\cdot)$ and $f_k(\cdot,\z_k^{t+1})$ also imply the following upper bounds
\begin{align}
& g_k(\z_k^{t+1}) \leq g_k(\x_k^t) + \langle \nabla g_k(\x_k^t) , 
\z_k^{t+1} - \x^t_{ k }\rangle \nonumber \\ 
& \qquad+\frac{L_k}{2}\sum_{j \in	\N_k}\norm{x_{kj}^t-z_j^{t+1}}^{ 2 } \label{lipcong3} \\
& f_k(\x_k^t,\z_k^{t+1}) \leq f_k(\z_k^{t+1},\z_k^{t+1})+\frac{L_k}{2}\sum_{j \in 
	\N_k}\norm{x_{kj}^t-z_j^{t+1}}^{ 2 } \nonumber \\
 & + \langle \nabla_{\x_k} f_k(\z_k^{t+1},\z_k^{t+1}) , \x_k^{t} - \z^{t+1}_{ k }\rangle \label{lipcong4}\\
	&\leq g_k(\x_k^{t}) + \langle \nabla g_k(\z_k^{t+1}) - g_k(\x_k^t), \x_k^{t} - \z^{t+1}_{ k }\rangle \nonumber \\ 
	&\qquad +L_k\sum_{j \in 
	\N_k}\norm{x_{kj}^t-z_j^{t+1}}^{ 2 } \label{fbnd}
\end{align}
Therefore, from \eqref{fbnd}, it follows that
\begin{subequations}\label{unboundf}
\begin{align}\label{ulbound}
& u_{ k }(\x_k^{t},\z_k^{t+1},\z_k^{t+1},\y_k^t)-\ell_{ k 
}(\x_k^{t},\z_k^{t+1},\y_k^t)\nonumber \\& =f_k(\x_k^t,\z_k^{t+1}) - g_k(\x_k^t)  \\
&  \leq  \langle \nabla 
	g_k({\z}_k^{t+1})-\nabla g_k( \x_k^{t}) , \x_k^{t}-{\z}_k^{t+1} \rangle 
\nonumber \\&\qquad	+L_k\sum_{j \in \N_k}\norm{x_{kj}^t-z_j^{t+1}}^{ 2 }\\
	&  \le 2L_k\sum_{j \in \N_k}\norm{x_{kj}^t-z_j^{t+1}}^{ 2 } 
	\label{ulboundf2}\\
	& \le 4L_k \sum_{j \in \N_k}\left(\norm{x_{kj}^{t}- x_{ kj }^{ t+1 }}^{ 
		2 }+\norm{x_{kj}^{t+1}-z_j^{t+1}}^{ 2 }\right)\label{ulboundf3} 
	\end{align}
\end{subequations}
where \eqref{ulboundf2} follows from (A4) and \eqref{ulboundf3} from the 
use of triangle inequality. Finally, it holds that
\begin{align}
&\ell_{ k }( \x_{ k }^{ t+1 },\z_k^{t+1},\y_k^{ t })-\ell_{ k 
}(\x_k^{t},\z_k^{t+1},\y_k^{ t })\\& 
\le u_k(\x_k^{t+1},\z_k^{t+1},\z_k^{t+1},\y_k^t) - u_k( 
\x_k^{t},\z_k^{t+1},\z_k^{t+1},\y_k^t) \nonumber \\ &+ u_k( \x_k^{t},\z_k^{t+1},\z_k^{t+1},\y^t)-\ell_k( 
\x_k^{t},\z_k^{t+1},\y^t)\\
&\leq -{ \frac { \rho _{ k }-9L_k }{ 2 } \sum_{j\in 
			\N_k}\|x_{kj}^{t+1}-x_{kj}^{t}\|^2 } +\frac { 4 L_k }{  \rho _{ k }^2 } 
	\sum_{j\in 
		\N_k}\|y_{kj}^{t+1}-y_{kj}^{t}\|^2\nonumber \\ &+\frac{L_kT_k}{2}\sum_{d=0}^{T_k-1}\sum_{j\in 
		\N_k}\|z_{j}^{t+1-d}-z_{j}^{t-d}\|^2 \label{ysub_maj}\\
	&= -{\left( \frac { \rho _{ k }-9L_k}{ 2 }-\frac{8L_k^3}{\rho_k^2}\right) \sum_{j\in 
			\N_k}\|x_{kj}^{t+1}-x_{kj}^{t}\|^2 } \nonumber \\ &
			\hspace{-5mm}+\left(\frac { 8 L^3_k(T_k+1)}{
		\rho _{ k }^2 }+\frac{L_kT_k}{2}\right)\sum_{d=0}^{T_k-1}\sum_{j\in 
		\N_k}\|z_{j}^{t+1-d}-z_{j}^{t-d}\|^2 \label{ysub2_maj}
	\end{align}
	where \eqref{ysub_maj} follows from using the update for $y_{kj}^{t+1}$ and 
	\eqref{ysub2_maj} from \eqref{ybnd_maj}. Finally, summing \eqref{ysub2_maj} over $k 
	= 1, 2, \ldots, K$, the required bound is obtained. Since the expression for the third summand in \eqref{lemz} remains the same,  the decrease in the Lagrangian over 
	consecutive time slots $t=t_0+1, \ldots, t_0+T$, is given by
	\begin{align} \label{convfint2}
 &	L( \{\x_{ k }^{ T+ t_0 }\},{\z}^{T+t_0},\{\y_k^{T+t_0}\}) -L( \{\x_{ k 
	}^{ t_0 }\},{\z}^{t_0}, \{\y_k^{t_0}\})\nonumber \\
	&\hspace{-5mm}\leq - \sum_{t=t_0}^{T+t_0}\sum_{k=1}^K { \beta_k 
		\sum_{j\in \N_k}\norm{x_{kj}^{i+1}-x_{kj}^{i}}^2 } + 
	\alpha_k\norm{z_{k}^{i+1}-z_{k}^{i}}^2 
	\end{align}
	where, $\alpha_k$ and $\beta_k$ are given in \eqref{betak2}
	yielding the desired result. Note that from Assumption (A5), the term on 
	the right-hand side of \eqref{convfint} is negative. Next, the boundedness of the Lagrangian follows as shown in Appendix \ref{{lem1}_proof}. Finally, it is possible to apply	Theorem \ref{thm1} to this case, yielding the desired result.

\section{Lipschitz continuity of objective function in \eqref{mds}}\label{lips}
Recall that $g_k(\{\x_j\}_{j \in \N_k}) = \sum_{j }w_{kj}(\delta_{kj} - d_{kj}(\x_k,\x_j))^2$, where the modified definition of $ d_{kj}(\x_k,\x_j) = \sqrt{\norm{\x_k-\x_j} + \epsilon}$ is utilized and for a pair of node $(k,j)$ weight $w_{kj}=1$ for  $j \in \N_k$ otherwise 0. The gradient of $g_k$ is given by $\nabla_{\x_l} g_k(\{x_j\}_{j \in \N_k})=$
\begin{align}\label{nablag}
\begin{cases} 2 w_{ kj } (\delta_{ kj } - d_{ kj }(\x_k,\x_j))\frac { (\x_{ k } -\x_{ j }) }{ d_{ kj }(\x_k,\x_j)} & l \neq k\\
-\sum_{j } 2 w_{ kj} (\delta_{ kj } - d_{ kj }(\x_k,\x_j))\frac { (\x_{ k } -\x_{ j}) }{ d_{ kj }(\x_k,\x_j)} & l = k
\end{cases}
\end{align}

In order to show that the expression in \eqref{nablag} is Lipschitz, it suffices to prove that each element of the Hessian is bounded. Denoting the coordinate vector of the $i$-th node as $\x_i := [x_i^1 ~ x_i^2]$ and defining the $2 \times 2$ matrix $[C(\x_k,\x_j)]_{pq} := 2w_{kj}(1-\delta_{kj}((x_k^p-x_j^p)(x_k^q-x_j^q)+\epsilon)d_{ kj }^{-3}(\x_k,\x_j))$ for $1 \leq p,q \leq 2$, it can be seen that the $(\ell,m)$-th block of the Hessian matrix is given by
\begin{align}\label{hess1}
\nabla^2_{\x_l,\x_m} g_k(\{x_j\}_{j \in \N_k}) =& \begin{cases}C(\x_k,\x_m)& m=l, l \neq k\\
-C(\x_k,\x_m),& m \neq l, m = k\\
{ \sum_{j\in \N_k' }}  C(\x_k,\x_j) & m=l,l=k \\
-C(\x_k,\x_l)  &  m\neq l,l=k\\
0 &m\neq l, l \neq k 
\end{cases}\nonumber 
\end{align}
For the present localization example, we have $\max_{k,j}\{\delta_{kj}\}\leq1$, $\max_{k}|\N_k|\leq N$, $w_{kj}\leq1$, $\| \x_k\|\leq 1\quad \forall \quad k=1,\ldots, K$, and $\| \frac{1}{d_{kj}(\x_k,\x_j)}\|\leq \frac{1}{\sqrt{\epsilon}}$. The application of triangle inequality therefore implies that $\norm{C(\x_k,\x_m)} \leq 8(2/\sqrt{\epsilon}+1)$ for all $1\leq  k,m \leq N$. Therefore it follows that $\nabla^2_{\x_l,\x_m} g_k(\{x_j\}_{j \in \N_k})$ is bounded, as claimed. 

\footnotesize
\bibliographystyle{IEEEtran}
\bibliography{IEEEabrv,reference}

\end{document}